\documentclass{amsart}

\usepackage{amsmath,amssymb, amscd}

\theoremstyle{plain}

\newtheorem{theorem}{Theorem}[section]
\newtheorem{proposition}[theorem]{Proposition}
\newtheorem{lemma}[theorem]{Lemma}
\newtheorem{corollary}[theorem]{Corollary}
\newtheorem{conjecture}[theorem]{Conjecture}

\theoremstyle{definition}

\newtheorem{definition}[theorem]{Definition}

\newtheorem{remark}[theorem]{Remark}

\newcommand{\ram}{\text{CR}}

\newcommand{\C}{\mathbf{C}}
\newcommand{\F}{\mathbf{F}}
\renewcommand{\H}{\mathcal H}
\newcommand{\I}{\mathcal{I}}
\renewcommand{\L}{\mathcal{L}}
\newcommand{\m}{\mathfrak m}
\newcommand{\M}{\mathcal M}
\renewcommand{\O}{\mathcal O}
\newcommand{\p}{\mathfrak p}

\renewcommand{\P}{\mathfrak P}

\newcommand{\Q}{\mathbf Q}
\newcommand{\T}{\mathbb T}
\newcommand{\X}{\mathcal X}
\newcommand{\Z}{\mathbf Z}

\newcommand{\GQ}{G_{\Q}}
\newcommand{\Ki}{K_{\infty}}
\newcommand{\Qinf}{\Q_{\infty}}

\newcommand{\Qlur}{\Q_{\ell}^{\text{ur}}}
\newcommand{\Qp}{\Q_{p}}
\newcommand{\Qbar}{\bar{\Q}}
\newcommand{\Qpbar}{\Qbar_{p}}
\newcommand{\Fpbar}{\overline{\F}_{p}}
\newcommand{\rhob}{\bar{\rho}}
\newcommand{\rhobar}{\rhob}
\newcommand{\vep}{\varepsilon}
\newcommand{\Zp}{\Z_{p}}

\newcommand{\Ifl}{\I_{f_{\ell}}}
\newcommand{\Jl}{J^{(\ell)}}

\newcommand{\inj}{\hookrightarrow}
\newcommand{\maps}{\longrightarrow}

\newcommand{\LL}{\mathfrak L}

\newcommand{\kie}{\Kinf,E}
\newcommand{\kif}{\Kinf,f}
\newcommand{\kig}{\Kinf,g}
\newcommand{\kih}{\Kinf,h}

\newcommand{\kifn}{\Kinf,f_{n}}

\newcommand{\lp}{L_{p}}
\newcommand{\Lp}{{\mathfrak L}_{p}}

\newcommand{\Kinf}{K_\infty}
\newcommand{\lf}{\lp(\kif)}
\renewcommand{\lg}{\lp(\kig)}

\newcommand{\lgn}{\lp(\Kinf,g)}
\newcommand{\lfpm}{\lp^{\pm}(\kif)}
\newcommand{\lfe}{\lp^\varepsilon(\kif)}
\newcommand{\lfpmn}{\lp^{\pm}(\kif)}
\newcommand{\lfp}{\lp^{+}(\kif)}
\newcommand{\lfm}{\lp^{-}(\kif)}
\renewcommand{\le}{\lp(\kie)}

\newcommand{\Lf}{\Lp(\kif)}

\newcommand{\Lfpm}{\Lp^{\pm}(\kif)}

\newcommand{\Lfpmn}{\Lp^\pm(\kif)}

\newcommand{\Le}{\Lp(\kie)}

\renewcommand{\sf}{\sel(\kif)}
\newcommand{\sg}{\sel(\kig)}
\newcommand{\sh}{\sel(\kih)}
\newcommand{\sfn}{\sel(\kifn)}

\newcommand{\sfpm}{\sel^\pm(\kif)}
\newcommand{\sfnpm}{\sel^\pm(\kifn)}
\newcommand{\se}{\sel_{p}(\kie)}

\newcommand{\Sf}{\Sel(\kif)}
\newcommand{\Sfpm}{\Sel^\pm(\kif)}
\newcommand{\Se}{\Sel_{p}(\kie)}

\newcommand{\bp}[1]{\bigl(#1\bigr)}
\newcommand{\mb}[1]{\mu\bp{#1}}

\newcommand{\Ts}{\T_0(N^+\!,N^-)}
\newcommand{\Tn}{\T_0(N)}
\newcommand{\Tv}[2]{\T_0(#1,#2)}

\newcommand{\cng}[3]{\eta_{#1}(#2,#3)}
\newcommand{\cf}{\eta_f(N)}
\newcommand{\xf}{\xi_f(N^+\!,N^-)}

\newcommand{\hXJr}{\hat{\X}_r(J)}
\newcommand{\hXAr}{\hat{\X}_r(A)}
\newcommand{\XJr}{{\X}_r(J)}

\newcommand{\hPJr}{\hat{\Phi}_r(J)}
\newcommand{\hPAr}{\hat{\Phi}_r(A)}
\newcommand{\PJr}{{\Phi}_r(J)}

\newcommand{\hT}{\hat{\T}}

\DeclareMathOperator{\alg}{alg}

\DeclareMathOperator{\Frob}{Frob}
\DeclareMathOperator{\Gal}{Gal}
\DeclareMathOperator{\GL}{GL}
\DeclareMathOperator{\Hom}{Hom}
\DeclareMathOperator{\ord}{ord}
\DeclareMathOperator{\Sel}{\mathfrak{Sel}}
\DeclareMathOperator{\sel}{Sel}
\DeclareMathOperator{\Ta}{Ta}
\DeclareMathOperator{\Ann}{Ann}
\DeclareMathOperator{\chr}{char}

\DeclareMathOperator{\Pic}{Pic}
\DeclareMathOperator{\coker}{coker}
\DeclareMathOperator{\un}{un}

\DeclareMathOperator{\fin}{fin}
\DeclareMathOperator{\res}{res}
\DeclareMathOperator{\corank}{corank}

\begin{document}

\title{On anticyclotomic $\mu$-invariants of modular forms}
\author{Robert Pollack and Tom Weston}
\address[Robert Pollack]{Department of Mathematics, Boston University, Boston, MA}
\address[Tom Weston]{Dept.\ of Mathematics, University of
Massachusetts, Amherst, MA}

\email[Robert Pollack]{rpollack@math.bu.edu}
\email[Tom Weston]{weston@math.umass.edu}

\thanks{Supported by NSF grants DMS-0439264 and DMS-0440708}

\maketitle

\section{Introduction}
\label{sec:intro}
Let $E/\Q$ be an elliptic curve of squarefree level $N$.
Fix a prime $p \geq 5$ of good reduction and an imaginary quadratic field $K$ of discriminant prime to $pN$.
Write $N=N^+N^-$ with $N^+$ divisible only by primes which are split in $K/\Q$ and $N^-$ divisible only by inert primes.  If $N^{-}$ has an even number of prime
divisors, then it is known by work of Cornut \cite{cornut}
and Vatsal \cite{Vmu} that
$E$ has infinitely many independent non-torsion points defined over ring class fields of $K$ of $p$-power conductor.
In terms of Iwasawa theory, this implies that $\se$, the $p$-adic Selmer group of $E$ over the anticyclotomic $\Zp$-extension $\Ki$ of $K$,
has positive rank over the Iwasawa algebra $\Lambda$; correspondingly,  the anticyclotomic $p$-adic $L$-function $\le$ vanishes. 

When $N^{-}$ has an odd number of prime divisors,
the situation is more analogous to the cyclotomic Iwasawa theory of $E$.
In this case, the signs
of the functional equations of twists of $E$ suggest that $E$ has finite
rank over $\Ki$: more precisely, one
expects that $\le$ is non-zero and that $\se$ is a
cotorsion $\Lambda$-module.
(When $E$ is $p$-supersingular one must replace these objects by their $\pm$-variants defined in \cite{DI,IP,Kobayashi}.  Our
discussion below continues to hold for these Selmer groups but for simplicity
we will focus in the
introduction on the ordinary case.)
Furthermore, the main conjecture of Iwasawa theory predicts that  $\le \cdot \Lambda$ equals the characteristic ideal of  $\se^{\vee}$, the Pontryagin dual of $\se$.

Many of these facts are now known in this setting: work of Vatsal
\cite{Vnonzero} establishes the non-vanishing of $\le$, while
under the additional hypothesis that $E$ is $p$-isolated,
Bertolini--Darmon \cite{BD}
have established the cotorsionness of $\se$.  In fact,
in this case \cite{BD} shows one divisibility of the main
conjecture: the characteristic ideal of $\se^\vee$ divides $\le$.

Remarkably, one can use the wealth of information provided by Heegner points in the complementary indefinite case to prove results in the definite case whose cyclotomic analogues remain unproven.  Specifically, in the cyclotomic case,  it is a long standing conjecture of Greenberg that the $\mu$-invariant of $E$ vanishes if the $p$-torsion of $E$ is irreducible.  The anticyclotomic analogue of this statement can be deduced from the work of Vatsal and Bertolini--Darmon.  Precisely, \cite{Vmu} establishes the vanishing of the analytic $\mu$-invariant; the divisibility of \cite{BD} then yields the vanishing of the algebraic $\mu$-invariant as well.  (See Theorem \ref{thm:mualg}.)

Vatsal's work also indicates an unexpected divergence from the cyclotomic setting: there are in fact two natural normalizations of the anticyclotomic $p$-adic $L$-function depending on whether one uses Gross' period of \cite{Gross} or Hida's canonical period \cite{Hida}. The $p$-adic $L$-function $\le$ discussed above corresponds to the first of these; we write $\Le$ for the second. There is in fact a corresponding choice of Selmer groups: the usual elliptic Selmer group $\se$ and the Selmer group $\Se$ in the sense of Greenberg.  (These Selmer groups differ only in the defining local conditions at primes dividing $N^-$: the former uses a locally trivial condition while
the latter uses only a locally unramified condition.)  

One goal of this paper is to illuminate the difference between these two choices.  Specifically, both algebraically and analytically they differ
only in the $\mu$-invariant.  Vatsal has in fact also given a precise formula for $\mb{\Le}$ in terms of congruence numbers:
$$
\mb{\Le} = \ord_{p}\left(\frac{\eta_{f}(N)}{\xf} \right).
$$
We define the quantities in this formula precisely in Section \ref{sec:analytic}.  For now, let us comment that  $\eta_{f}(N)$ measures congruences between the newform $f$ corresponding to $E$ and other eigenforms
of weight two and level $N$. The term $\xf$
is closely related to congruences with such eigenforms that are also new at all primes dividing $N^-$.
(We remark that this formula differs slightly from the formula stated in \cite{Vmu}.  We will elaborate on this difference in Section \ref{sec:analytic}.)

In this paper we obtain a very different looking formula for $\mb{\Se}$.  We state this formula here in general for weight two modular forms.
Let $f$ be a newform of weight two and squarefree level $N=N^+N^-$ such that the number of prime divisors of $N^-$ is odd.  
Throughout the paper we will be imposing the following hypotheses on $N^-$ and $\rhob_f$, the residual representation attached to $f$:
\vspace{.2cm}
\begin{quote}
A continuous Galois representation $\rhob : \GQ \to \GL_{2}(\Fpbar)$ and a squarefree product $N^-$ of an odd number of primes, each inert in $K/\Q$, including all 
such primes at which $\rhob$ is ramified, satisfies {\it hypothesis $\ram$}~if:
\vspace{.2cm}
\begin{enumerate}
\item $\rhob$ is surjective;
\item \label{hyp:free}
if $q \mid N^-$ and $q \equiv \pm 1 \pmod{p}$, then $\rhob$ is ramified at $q$.
\end{enumerate}
\end{quote}
\vspace{.2cm}
In the below theorem, the Tamagawa exponent $t_{f}(q)$ is a purely local invariant which for an elliptic curve is simply the $p$-adic valuation of the Tamagawa factor at $q$.  See Definition~\ref{def:tam} for a precise description of this quantity in general.

\begin{theorem}  
Assume that $(\rhob_f,N^{-})$ satisfies $\ram$.
\label{thm:main}
\begin{enumerate}
\item If $f$ is $p$-ordinary, then
$$
\mb{\sf}=0 \text{~and~}
\mb{\Sf}= \sum_{q \mid N^-} t_{f}(q).$$
\item
If $f$ is $p$-supersingular, $a_{p}=0$, $p$ is split in $K/\Q$ and each prime above $p$ is totally ramified in $\Ki/K$, then
$$\mb{\sfpm}=0  \text{~and~}
\mb{\Sfpm} = \sum_{q \mid N^-} t_{f}(q).
$$
\end{enumerate}
\end{theorem}

We remark that these results are strikingly similar to results of
Finis \cite{Finis} on analytic $\mu$-invariants of anticyclotomic Hecke-characters.

Note that the formula of Vatsal gives the analytic $\mu$-invariant as a difference of global terms while the algebraic formulae above are purely local.   The main conjecture, however, predicts that these formulae should agree.  Using results of Ribet--Takahashi \cite{RT,T} on degrees of modular parameterizations arising from Shimura curves, we establish that 
\begin{equation}
\label{eqn:intro}
\ord_{p}\left(\frac{\eta_{f}(N)}{\xf} \right)
= \sum_{q \mid N^-} t_{f}(q)
\end{equation}
and thus deduce the $\mu$-part of the main conjecture in both the ordinary and supersingular case.

Another goal of this paper is to weaken the hypotheses of the results of Bertolini--Darmon \cite{BD} on the cotorsionness of $\sf$ and on the divisibility of $\lf$ by the characteristic power series of $\sf^\vee$.  In \cite{BD}, it is assumed that $f$ has Fourier coefficients in $\Zp$ and that $f$ is not congruent to any eigenform of level $N$ that is old at any prime dividing $N^-$.  These hypotheses are not stable under congruences and so are unfavorable for studying congruence questions in the spirit of \cite{GV,EPW}.

To remove the first assumption on the Fourier coefficients more care is needed in studying the Galois representations that arise; this is dealt with in Proposition \ref{prop:gal}.  The second assumption is used in two serious ways in \cite{BD}.  It is used to (trivially) deduce the freeness of a certain character group attached to a Shimura variety.  This freeness is used to carry out mod $p^n$ level-raising to produce a mod $p^n$ modular form congruent to $f$ with certain desirable properties.  It is then used again to lift this mod $p^n$ modular form to a true modular form.  

We address the character group via
hypothesis \ram~ for $(\rhob_f,N^-)$ (which is weaker than the hypotheses of \cite{BD}), showing that it is enough to force the freeness of the character group
(see Theorem \ref{thm:free}).  Thus, the mod $p^n$ level-raising arguments can still be made to work.  As for the second issue,
it is not possible in general to lift mod $p^n$ modular forms to true modular forms.  We circumvent this problem by working directly with mod $p^n$ modular forms, their Selmer groups and their $p$-adic $L$-functions.
One then verifies that the arguments of \cite{BD} go through in this more general setting.

We close this introduction by proposing a formula on congruences numbers that is purely a statement about modular forms, but arises naturally from the study of anticyclotomic $\mu$-invariants especially equation (\ref{eqn:intro}).  Namely,
if $N=a\ell b$ is a factorization of the level of $f$ with $\ell$ a prime, we conjecture that
\begin{equation} \label{eq:pf}
\vspace{.1cm}
\ord_p \bp{ \eta_{f}(a\ell,b)} = t_{f}(\ell) + \ord_p \bp{\eta_{f}(a,\ell b)}.
\end{equation}
Here, for a factorization $N=N_{1}N_{2}$, the quantity
$\eta_f(N_{1},N_{2})$ measures congruences between $f$ and forms of level $N$ that are new at all primes dividing $N_{2}$.
(See Section \ref{sec:pf} for a precise statement.)

This conjecture immediately implies level-lowering in the sense
of \cite{Ribet}; it perhaps should be regarded as a quantitative version of level-lowering, much as Wiles' numerical
criterion \cite{Wiles} is a quantitative version of level-raising.
The formula of Ribet--Takahashi referred to above is an analogue of this formula in terms of degrees of modular parameterizations arising from Shimura curves.  Similar formulae appear in Khare's work \cite{Khare} on establishing isomorphisms between deformation rings and Hecke rings via level-lowering.
We prove (\ref{eq:pf}) in Section \ref{sec:pf} assuming $\ram$.
Not coincidently,  this hypothesis puts us in the case in which level-lowering can be established by Mazur's principle.

The structure of the paper is as follows.
In Section 2, we define the two normalizations of $p$-adic $L$-functions and recall the results of Vatsal.  In Section 3, we define our two normalizations of Selmer groups and compare them.  In Section 4, we generalize the results of \cite{BD} as described above.  In Section 5, we combine the results of the previous two sections to produce a local formula for the algebraic $\mu$-invariants.  In Section 6, we prove the $\mu$-part of the main conjecture and discuss quantitative level-lowering.  Finally, in section 7, we discuss some anticyclotomic analogues of the congruence results of \cite{GV,EPW}.  We also include an appendix giving a general criterion for surjectivity of
global-to-local maps in Iwasawa theory.

\vspace{.2cm}
\noindent
{\it Acknowledgments:}  
We would like to thank Henri Darmon, Matthew Emerton, Ralph Greenberg, Farshid Hajir, and Christian Maire for their help with various aspects of this paper.

\section*{Notation}
Fix an odd prime $p$ and embeddings $\Qbar \inj \Qpbar$ and $\Qbar \inj \C$.  Let $K/\Q$
be an imaginary quadratic field with discriminant $D$ prime to $p$.
Let $\Kinf$ denote the anticyclotomic
$\Zp$-extension of $K$.  Thus $\Gamma := \Gal(\Kinf/K)$ is non-canonically isomorphic
to the additive group $\Zp$ and the non-trivial element of
$\Gal(K/\Q)$ acts on $\Gamma$ by inversion.  We write $K_{n}$
for the unique subfield of $\Kinf$ such that $\Gal(K_{n}/K) \cong
\Z/p^{n}\Z$.  If $N$ is an integer relatively prime to $D$, we write
$N^{+}$ (resp.\ $N^{-}$) for the largest divisor of $N$ divisible only
by primes split (resp.\ inert) in $K/\Q$.

Let $f = \sum a_n q^n$ denote 
a normalized newform of weight two, squarefree level $N=N^+N^-$
prime to $pD$, and trivial nebentypus.  
We assume throughout this paper that $N^-$ has an {\it odd} number of prime
factors.
We regard $f$ as a $p$-adic modular form via our fixed
embedding $\Qbar \inj \Qpbar$; let $\O_{0}$ denote the $\Zp$-subalgebra of
$\Qpbar$ generated by the images of the Fourier coefficients of $f$ and let
$\O$ denote the integral closure of $\O_{0}$ in its fraction field $F$.  
We write $\p$ for
the maximal ideal of $\O$ and for $n \geq 1$ set $\p_{n} := \p^{n} \cap
\O_{0}$.  Let $k_{0} = \O_{0}/\p_{1}$ and $k = \O/\p$ denote the residue
fields.
Let $\Lambda_{0} := \O_{0}[[\Gamma]]$ and $\Lambda := \O[[\Gamma]]$
denote the Iwasawa algebra over $\O_{0}$ and $\O$ respectively.
 
\section{$p$-adic $L$-functions and analytic $\mu$-invariants}
\label{sec:analytic}

\subsection{The complex period $\Omega$}

When $f$ is $p$-ordinary, in the sense that $a_{p}$ is a $p$-adic unit, there is an anticyclotomic $p$-adic $L$-function
$$\lf \in \Lambda$$
interpolating the algebraic special values of the $L$-series of anticyclotomic twists of $f$ over $K$.  In particular, for a character $\chi$ of $\Gamma$ of order $p^{n}$ we have (up to $p$-adic units)
\begin{eqnarray}
\label{eqn:interpolate}
\chi\bp{\lf} = 
\frac{1}{\alpha^{2n}} \cdot \frac{L(f,\chi,1)}{\Omega} \cdot C_\chi,
\end{eqnarray}
where $\alpha$ is the unit root of $x^{2}-a_{p}x + p$,
$C_\chi = \sqrt{D}p^n$, and $\Omega := \Omega_{f,K}$ is a certain complex period that 
depends upon $f$ and $K$ as in \cite{BD96,BDCD,BD,Vmu}. 
We recall now the definition of $\Omega$.

Fix a factorization $N=N_{1}N_{2}$ and let $S_2(N_1,N_2)$ denote the space of cusp forms on $\Gamma_0(N)$ that are new at all primes dividing $N_2$.  Let $\Tv{N_{1}}{N_{2}}$ denote the $p$-adic completion of the Hecke algebra that acts faithfully on $S_2(N_1,N_2)$.  We simply write $\Tn$ for $\Tv{N}{1}$.  Our fixed newform $f$ gives rise to a homomorphism
\begin{equation} \label{eq:hecke}
\pi_f : \Tv{N_{1}}{N_{2}} \maps \O_0
\end{equation}
sending $T_\ell$ to $a_\ell$ for every prime $\ell \nmid N$ and $U_q$ to $a_q$ for each prime $q \mid N$. 

Let $X_{N^+\!,N^-}$ denote the Shimura curve of level $N^+$ attached to the definite quaternion algebra ramified at the primes dividing $N^-$.  If $\M = \Pic(X_{N^+\!,N^-}) \otimes \Zp$, then $\M$ has a natural 
faithful action of $\Ts$.
In the construction of the $p$-adic $L$-function $\lf$, one chooses a linear
map
$$\psi_f : \M \to \O$$
that is $\Ts$-equivariant where $\O$ is viewed as a $\Ts$-module via $\pi_f$.   As $\M \otimes \Qp$ is a free $\Ts \otimes \Qp$-module, this map is uniquely determined up to multiplication by an element of $\O$.  By scaling by a constant of $\O$, we can and do insist that $1$ be in the image of $\psi_f$.  This normalization determines $\psi_f$ up to a $p$-adic unit.

We now explicitly construct such a map.   Let $\M^f$ denote the submodule of $\M \otimes \O$ on which $\Ts$ acts via $\pi_f$.  Then $\M^f$ is a free $\O$-module of rank 1; let $g_f$ denote a generator of this module.  The Hecke-module $\M$ is equipped with an intersection pairing $\langle \cdot , \cdot \rangle : \M \times \M \to \Zp$ under which the action of $\Ts$ is adjoint.

\begin{lemma}
\label{lemma:unit}
There is some $m \in \M$ such that $ \langle m, g_f \rangle$ is a unit.
\end{lemma}

We prove this lemma by relating $\M$ to a character group arising from a Shimura curve attached to an indefinite quaternion algebra and then invoking results of \cite{T}.  As such character groups will be explored in detail in Section \ref{sec:char}, we postpone a proof until then.

Assuming this lemma, we may
take the map $\psi_f$ to be defined by 
$$
\psi_f(x) = \langle x , g_f \rangle.
$$
Under this choice of normalization, we now specify the period $\Omega$ precisely.  Set 
$$
\xf = \langle g_f , g_f \rangle.
$$
As $g_f$ is only defined up to a $p$-adic unit, we can choose $g_f$ so that $\xf$ is in $K$ and thus we may view it as an element of $\Qpbar$ or $\C$ via our fixed embeddings.  
  
\begin{lemma}
\label{lemma:period}
The period $\Omega$ in (\ref{eqn:interpolate}) can be taken to be
$$
\Omega = \frac{ (f,f) }{\xf}.
$$
Here $(f,f)$ denotes the Petersson inner product of $f$ with itself.
\end{lemma}

\begin{proof}
Let $e_f$ be the idempotent of $\Ts \otimes \Qp$ attached to $f$.  Let $P \in \M$ (a Heegner point) and $G_n$ (a Galois group) be as in \cite[Lemma 2.5]{Vmu}. Gross' special value formula implies that for $\chi$ a primitive character of $G_n$, we have
$$
\left\langle 
\sum_{\sigma \in G_n} \chi(\sigma) \cdot e_f  P^\sigma ,
\sum_{\sigma \in G_n} \chi(\sigma) \cdot e_f P^\sigma  
\right\rangle
=
\frac{L(f,\chi,1)}{(f,f)} \cdot C_\chi \cdot u^2
$$
where $u$ is half the number of units of $\O_K^\times$.
Expanding the left hand side gives
$$
\sum_{\sigma,\tau} \chi(\sigma \tau^{-1}) 
\langle  e_f P^\sigma  , e_f P^\tau  \rangle
=
\frac{L(f,\chi,1)}{(f,f)} \cdot C_\chi \cdot u^2.
$$
Write $e_f P^\sigma = c_f(\sigma) g_f$ with $c_f(\sigma)$ in the field of fractions of $\O$.  (Note that $c_f(\sigma)$ need not be integral as $e_f P^\sigma$ is an element of $\M \otimes \Qp$ and not necessarily of $\M$.)  Then
$$
\psi_f(e_f P^\sigma) = \langle e_f P^\sigma, g_f \rangle = c_f(\sigma) \langle g_f, g_f \rangle
$$
and thus
$$
\sum_{\sigma,\tau} \chi(\sigma \tau^{-1}) \psi(e_f P^\sigma) \psi(e_f P^\tau) = 
\frac{L(f,\chi,1)}{(f,f)} \cdot C_\chi \cdot u^2 \cdot \langle g_f, g_f \rangle.
$$
Rearranging yields
$$
\left| \sum_\sigma \chi(\sigma) \psi(e_f P^\sigma) \right|^2  =
\frac{L(f,\chi,1)}{(f,f)} \cdot C_\chi \cdot u^2 \cdot \langle g_f, g_f \rangle
$$
and since
$$
\psi(e_f P^\sigma) = \langle e_f P^\sigma, g_f \rangle = \langle P^\sigma, e_f g_f \rangle = \langle P^\sigma,g_f \rangle = \psi(P^\sigma)
$$
we have
$$
\left| \sum_\sigma \chi(\sigma) \psi(P^\sigma) \right|^2  =
\frac{L(f,\chi,1)}{(f,f)} \cdot C_\chi \cdot u^2 \cdot \langle g_f, g_f \rangle.
$$
By the definition of the $p$-adic $L$-function, it then follows that 
we may take
$$	
\Omega = \frac{(f,f)}{\xf}.
$$
\end{proof}

\subsection{The canonical period}

There is a second natural choice of complex period in this context, namely Hida's canonical period (cf.\ \cite{Hida}, \cite[pg.\ 10]{Vmu}).  
In order to state the definition of this period, we recall the notion of a congruence number.

For later use we proceed in somewhat more generality than is immediately necessary.  For a factorization $N=N_{1}N_{2}$,
we define the congruence number $\cng{f}{N_{1}}{N_{2}}$ to be any generator of the $\O$-ideal
$$\pi_f\bp{\Ann_{\Tv{N_{1}}{N_{2}}}(\ker \pi_f)} \cdot \O;$$
here $\pi_f$ is as in (\ref{eq:hecke}).
The congruence number $\eta_{f}(N_{1},N_{2})$ 
is a unit if and only if $f$ does not admit any non-trivial congruences to eigenforms for $\Gamma_0(N)$ that are new at every prime dividing $N_{2}$. We simply write $\cf$ for $\cng{f}{N}{1}$.  (As $\cng{f}{N_{1}}{N_{2}}$ is only defined up to a $p$-adic unit, we can choose it to be in $K$ and thus view it in either $\Qpbar$ or $\C$.)

The {\it canonical period} of $f$ is defined as
$$
\Omega_f := \frac{ ( f,f )}{\cf}
$$
with $( f,f )$ as before the Petersson inner-product of $f$ with itself.  
This is a natural period to consider 
from the point of view of congruences; see, for example,
\cite{Vcong}.  Note also that it is independent of the imaginary
quadratic field $K$.

We denote the anticyclotomic $p$-adic $L$-function of $f$
relative to the canonical period $\Omega_f$ by $\Lf$; thus (up to $p$-adic units)
$$
\chi\bigl(\Lf\bigr) = \frac{1}{\alpha^{2n}} \cdot \frac{L(f,\chi,1)}{\Omega_f} \cdot C_\chi,
$$
with notation as before. By Lemma \ref{lemma:period}, we can choose $\Omega_f$ so that
\begin{eqnarray}
\label{eqn:lfLf}
\Lf = \lf \cdot \frac{\eta_f(N)}{\xf}.
\end{eqnarray}

\subsection{Analytic $\mu$-invariants}

It is known by work of Vatsal \cite{Vnonzero} that the $p$-adic $L$-functions $\lf$ and $\Lf$ are non-zero.  In fact, the results of \cite{Vmu} give the precise value of their Iwasawa $\mu$-invariants.  (We normalize our $\mu$-invariants so that $\mu(Q)$ for $Q \in \Lambda$ is the largest exponent $c$ such that
$Q \in \p^{c} \cdot \Lambda$.)

\begin{theorem}[Vatsal]
\label{thm:muanal}
Assume that the residual Galois representation $\rhobar_f$ 
attached to $f$ is irreducible.  Then:
\begin{enumerate}
\item $\mb{\lf} = 0$;
\item $\displaystyle \mb{\Lf} 
= \ord_\p \left(\frac{\cf}{ \xf} \right)$.
\end{enumerate}
\end{theorem} 

\begin{remark}
In \cite{Vmu}, the denominator $\xf$ is mistakenly replaced with the congruence number $\eta_f(N^+\!,N^-)$.  These two quantities are equal if $\M$ is a free $\Ts$-module (see Theorem \ref{thm:free}) and are probably always equal, but this is currently not known.  The relation between these two quantities will be further explored in Section \ref{sec:pf}.
\end{remark}

\begin{proof}
The first part follows from the discussion in \cite[Section 4.6]{Vmu}; see, in particular, \cite[Proposition 4.7]{Vmu}.  (Note that $\nu=0$ in the notation of \cite{Vmu} because we are assuming $\rhobar_f$ is irreducible.)  The second part follows immediately from (\ref{eqn:lfLf}).
\end{proof}

\subsection{Supersingular case}

We now turn to the case where $f$ is $p$-supersingular. Under the assumption that $a_{p} = 0$, two $p$-adic $L$-functions
$$\lfp, \lfm \in \O[[\Gamma]]$$
are constructed in \cite{DI} (in an analogous way to the cyclotomic construction of \cite{P}). These $p$-adic $L$-functions are defined with respect to the period $\Omega$ above; as before we denote by $\Lfpm$ the corresponding $p$-adic $L$-functions normalized with respect to the canonical period $\Omega_{f}$. If we only have that $a_{p} \equiv 0 \pmod{\p_{n}}$,  then one may still construct $p$-adic $L$-functions
$$\lfpmn, \Lfpmn \in \O/\p^n[[\Gamma]].$$
In the case when $a_p=0$, these $L$-functions are simply the mod $\p_{n}$ reductions of $\lfpm$ and $\Lfpm$.

The results of \cite{Vmu} extend easily to the supersingular case.

\begin{theorem} \label{thm:muanalss}
Let $f$ be as above and assume that $a_p = 0$.  Then:
\begin{enumerate}
\item $\mb{\lfpm} = 0$;
\item $\displaystyle \mb{\Lfpm} = \ord_{\p} 
\left(\frac{\cf}{\xf }\right)$.
\end{enumerate}
\end{theorem} 

\begin{proof}
In \cite{DI}, the $p$-adic $L$-functions $\lfpm$ are constructed via a sequence
$\{ \L_n \}_{n \geq 1} $ (with  $\L_n \in \O[\Gal(K_n/K)]$) 
satisfying $\pi^{n}_{n-1} \L_n = -\xi_{n-1} \L_{n-2}$; here
$$\pi^n_{n-1} : \O[\Gal(K_{n}/K)] \to \O[\Gal(K_{n-1}/K)]$$
is the natural map  and $\xi_{n-1} = \sum_{\sigma \in \Gal(K_{n-1}/K_{n-2})} \sigma$.

It follows that $\L_n$ is divisible by $\omega_{n}^{\vep}$ where
$$
\omega_n^+ = \prod_{\genfrac{}{}{0pt}{2}{1 \leq k \leq n}{k \text{~even}}}
\xi_k
\hspace{.5cm} \text{and} \hspace{.5cm}
\omega_n^- = \prod_{\genfrac{}{}{0pt}{2}{1 \leq k \leq n}{k \text{~odd}}} \xi_k,
$$
and $\vep$ equals the sign of $(-1)^{n-1}$.
For a fixed parity of $n$, factoring out these extra zeroes then 
produces the norm compatible sequence that yields $\lfe$.

The arguments of  \cite[section 5.9]{Vmu} in the ordinary case 
(which make use of \cite[Propositions 4.7 and 5.6]{Vmu})
generalize immediately to show that 
$\mu(\L_n)=0$ for $n$ large enough.  
Since $\mu(\omega^\pm_n)=0$ for all $n$, we deduce that 
$\mb{\lfpm}=0$
as desired.  The second part follows from (\ref{eqn:lfLf}).
\end{proof}

\section{Selmer groups}
\label{sec:alg}

We continue with the notation of the previous section.
Let $V_{f}$ denote the $p$-adic Galois representation associated to $f$: it is
a two-dimensional $F$-vector space endowed with a continuous action of
$\GQ$.  Fix an $\O$-stable lattice $T_{f} \subseteq V_{f}$ and set
$A_{f} = V_{f}/T_{f}$.  
We assume throughout that the residual
Galois representation 
$T_{f}/\p T_{f} \cong A_{f}[\p]$ is absolutely irreducible,
in which case $T_{f}$ is uniquely determined up to scaling.

\subsection{$p$-adic Selmer groups}
As with $p$-adic $L$-functions, there are two natural notions
of $p$-adic Selmer groups of $f$ over the anticyclotomic $\Zp$-extension
$\Ki$ of $K$: at places different from $p$, one can require that the
cocycles be either
locally trivial  or locally unramified.

Assume initially 
that $f$ is a $p$-ordinary modular form so that $V_f$ is a $p$-ordinary
Galois representation.
Define the {\it minimal Selmer group} $\sf$ as the kernel of
$$
H^1(\Ki,A_f) \maps \prod_{w \nmid p}  H^1(K_{\infty,w},A_f) \times
\prod_{w \mid p}  \frac{H^1(K_{\infty,w},A_f)}{H^1_{\ord}(K_{\infty,w},A_f)}
$$
and the {\it Greenberg Selmer group} $\Sf$ as the kernel of
$$
H^1(\Ki,A_f) \maps \prod_{w \nmid p}  H^1(I_{\infty,w},A_f) \times
\prod_{w \mid p}  \frac{H^1(K_{\infty,w},A_f)}{H^1_{\ord}(K_{\infty,w},A_f)};
$$
here $w$ runs over all places of $\Ki$, $I_{\infty,w}$ denotes the inertia group at $w$ and $H^1_{\ord}(K_{\infty,w},A_f)$ is the standard ordinary condition of \cite[p.\ 98]{Gord}.

\begin{remark}
If $w$ is not split infinitely in $\Ki$, then $G_{K_{\infty,w}}/I_{\infty,w}$ has profinite degree prime to $p$.  In particular, the map $H^1(K_{\infty,w},A_f)$ to $H^1(I_{\infty,w},A_f)$ is injective and the local condition at $w$ agree for both of these Selmer groups.

Since there are no primes which split infinitely in a cyclotomic $\Zp$-extension, the minimal and Greenberg Selmer groups coincide in this case.  In the
anticyclotomic case, the existence of primes which do split infinitely (namely,
those which are inert in $K/\Q$) can cause the Greenberg Selmer group to be
strictly larger than the minimal Selmer group.
\end{remark}

To compare these two Selmer groups, we give a more explicit description of the local conditions defining them.  
For $\ell$ a prime and $m \leq \infty$, let
$\sigma_{m,\ell}$ denote the set of places of $K_{m}$ lying over $\ell$ and set
$$\H_{\ell} = \underset{m}{\varinjlim} \prod_{w \in \sigma_{m,\ell}} H^{1}(K_{m,w},A_{f}) \hspace{.2cm} \text{~and~} \hspace{.2cm} \H_p = \prod_{w \in \sigma_{\infty,p}}  \frac{H^1(K_{\infty,w},A_f)}{H^1_{\ord}(K_{\infty,w},A_f)}$$
for $\ell \neq p$.

\begin{lemma} \label{l:hv}
Let $\ell \neq p$ be prime.
If $\ell$ is split in $K/\Q$, then $\sigma_{\infty,\ell}$ is finite and
$$\H_{\ell} = \prod_{w \in \sigma_{\infty,\ell}} H^{1}(K_{\infty,w},A_{f}).$$
If $\ell$ is inert or ramified in $K/\Q$, then $\ell$ splits completely in $\Ki$ and
$$\H_{\ell} = H^{1}(K_{\ell},A_{f}) \otimes \Lambda^{\vee}$$
where $\Lambda^{\vee} = \Hom_{\O}(\Lambda,F/\O)$.
\end{lemma}
\begin{proof}
The descriptions of $\sigma_{\infty,\ell}$ follow easily from class field
theory.  The formula for $\H_{\ell}$ in the split case follows immediately.
In the inert case, we have $K_{m,w} = K_{\ell}$ for all $w \in
\sigma_{m,\ell}$, so that
$$\prod_{w \in \sigma_{m,\ell}} H^{1}(K_{m,w},A_{f}) \cong
H^{1}(K_{\ell},A_{f}) \otimes \O[\Gal(K_{m}/K)].$$
Taking the limit over $m$ yields the desired description of $\H_{\ell}$ in this case.
\end{proof}

We now compute the difference between the defining local conditions of the
minimal and Greenberg Selmer groups.
For a prime $\ell$ that is inert in $K$, let 
$\H_{\ell}^{\un} \subseteq \H_{\ell}$ denote the
set of unramified cocycles; that is, the kernel of the map
$$H^{1}(K_{\ell},A_{f}) \otimes \Lambda^{\vee}
\to H^{1}(I_{\ell},A_{f}) \otimes \Lambda^{\vee}$$
with $I_{\ell} \subseteq G_{K_\ell}$ the inertia group at $\ell$.
This subgroup can be computed quite explicitly.  We first make a definition.

\begin{definition}
\label{def:tam}
Let $\ell$ be a prime number.  We define the
{\it Tamagawa exponent} $t_{f}(\ell)$ as follows.
If $A_{f}$ is unramified at $\ell$, we set $t_{f}(\ell)=0$.  If $A_{f}$
is ramified at $\ell$, we let $t_{f}(\ell)$ denote the largest exponent
$t \geq 0$ such that $A_{f}[\p^{t}]$ is unramified at $\ell$.
\end{definition}

\begin{lemma} \label{lemma:tam}
Let $\ell$ be a prime that is inert in $K$.  Then
$$\H_{\ell}^{\un} \cong  \O/\p^{t_f(\ell)} \otimes \Lambda^{\vee}.$$
In particular, $\mu(\H_{\ell}^{\un}) = t_{f}(\ell)$ and
$\lambda(\H_{\ell}^{\un})=0$.
\end{lemma}

\begin{proof}
From the inflation--restriction sequence, we have
$$\H_{\ell}^{\un} \cong H^{1}(k_\ell,A_{f}^{I_{\ell}}) \otimes \Lambda^{\vee}$$ 
with $k_\ell = \F_{\ell^2}$ the residue field at $\ell$.
Since the absolute Galois group of a finite field is pro-cyclic, we have
$$H^{1}(k_\ell,A_{f}^{I_{\ell}}) \cong A_{f}^{I_{\ell}}/(\Frob_{\ell}-1)A_{f}^{I_{\ell}}$$
with $\Frob_{\ell}$ a Frobenius element at $\ell$.  If $A_{f}$ is unramified
at $\ell$, then $A_{f}^{I_{\ell}} = A_{f}$ is divisible and has no trivial
Frobenius
eigenvalues, from which it follows that $H^{1}(k_{v},A_{f}^{I_{v}})$
vanishes.  This proves the lemma in the unramified case.

If $A_{f}$ is ramified at $\ell$, then
$$A_{f} \cong \left( \begin{array}{cc} \vep\chi & * \\ 0 & \chi \end{array}
\right)$$
as a $G_{\Q_\ell}$-module;
here $\vep$ is the cyclotomic character and $\chi$ is an unramified
quadratic character.
Since $\ell$ is inert in $K/\Q$, it follows that as a $G_{K_\ell}$-module we have
$$A_{f} \cong \left( \begin{array}{cc} \vep & * \\ 0 & 1 \end{array}
\right).$$
By the definition
of the Tamagawa exponent, it follows that
$$A_{f}^{I_{\ell}} = F/\O(\vep) \oplus \O/\p^{t_f(\ell)}.$$
The Frobenius coinvariants of this module simply equal $\O/\p^{t_f(\ell)}$,
as desired.
\end{proof}

 It follows from the definitions above
that there is an exact sequence
$$0 \to \sf \to \Sf \to \prod_{\ell \mid N^{-}} \H_{\ell}^{\un}.$$
In particular, it follows from Lemma~\ref{lemma:tam}
that the minimal Selmer group is $\Lambda$-cotorsion if
and only if the Greenberg Selmer group 
is $\Lambda$-cotorsion, and that when this is the
case, they have equal $\lambda$-invariants.
We will prove later that these 
sequences are in fact exact on the right as well, which allows
for a comparison of their $\mu$-invariants.

If $f$ is a  $p$-supersingular modular form, then 
the above discussion goes through if for each place $w$ of $\Ki$ dividing $p$ one replaces $H^1_{\ord}(K_{\infty,w},A_f)$
with the plus/minus local condition of \cite{Kobayashi,IP}.  Such a condition is
presently only defined under the assumptions that $a_{p}=0$, 
$p$ is split in $K$  and each prime above $p$ is totally ramified in $\Ki/K$.  We will make these assumptions from now on whenever
dealing with the $p$-supersingular case.
In particular, the above discussion yields Selmer groups $\sfpm$ and $\Sfpm$ together with exact sequences
$$0 \to \sfpm \to \Sfpm \to \prod_{\ell \mid N^{-}}
\H_{\ell}^{\un}$$
for each choice of sign.  The analysis of $\lambda$ and $\mu$-invariants applies equally well to this setting.

\subsection{Residual Selmer groups}
\label{sec:modpsel}

Assuming that $(\rhob_f,N^-)$ satisfies hypothesis
$\ram$ of the introduction, we define a {\it residual} Selmer group
$\sfn$ of $A_{f}[\p^{n}]$, depending on the Galois module
$A_{f}[\p^{n}]$ and $N^{-}$ but not on $f$ itself.
 
\begin{lemma}
\label{lemma:line}
Assume that $\ram$~ holds for ($\rhob_f,N^{-})$.  Then for each prime $\ell \mid
N^{-}$ there is a unique free rank one $\O/\p^{n}$-direct summand
$A^{(\ell)}_{f,n}$ of $A_{f,n}$ on which $G_{\Q_\ell}$ acts by either the
cyclotomic character or its negative.
 \end{lemma}
 
\begin{proof}
This is immediate from $\ram$~ for $n=1$; the general case follows easily
by induction.
\end{proof}

Define $\H_{\ell,n}$ as in the last section by replacing $A_f$ with $A_{f,n}$.  We also define
$$\H_{p,n} = \prod_{w \in \sigma_{\infty,p}}  \frac{H^1(K_{\infty,w},A_f[\p^{n}])}{H^1_{\ord}(K_{\infty,w},A_f[\p^{n}])}.$$
If $\ell$ is inert in $K$, let $\H^{\ord}_{\ell,n}$ denote 
the image of $H^1(K_\ell,A_{f,n}^{(\ell)}) \otimes \Lambda^\vee$ in $\H_{\ell,n}$. 

Fix a set $\Sigma$ of places of $\Q$ containing all primes dividing $Np$ and all archimedean places, but no primes that ramify in $K/\Q$.
Let $\Sigma^{-}$  (resp.\ $\Sigma^+$) denote the subset of primes of $\Sigma - \{p\}$ which are inert (resp.\ split) in $K$.  We define the residual Selmer group $\sfn$ as the kernel of
$$
H^1(K_\Sigma/\Ki,A_{f,n}) \to \H_{p,n} \times \prod_{\ell \in \Sigma^+} \H_{\ell,n} \times \prod_{\ell \in \Sigma^-} \H_{\ell,n} / \H^{\ord}_{\ell,n}.
$$
We note that these are the Selmer groups that are considered in \cite{BD}. 
If $f$ is $p$-supersingular, $a_{p}=0$, $p$ is split in $K/\Q$  and each prime above $p$ is totally ramified in $\Ki/K$, we may define $\sfnpm$ analogously as in \cite{DI}.

\begin{proposition} \label{prop:limit}
Assume that $(\rhob_f,N^-)$ satisfies $\ram$.
If $f$ is $p$-ordinary, then
$$
\sf = \underset{n}{\varinjlim} \, \sfn.
$$ 
If $f$ is $p$-supersingular, $a_{p}=0$, $p$ is split in $K/\Q$  and each prime above $p$ is totally ramified in $\Ki/K$, then 
$$
\sfpm = \underset{n}{\varinjlim} \, \sfnpm.
$$ 
\end{proposition}
 
\begin{proof}
We treat only the case when $f$ is $p$-ordinary as the proof in the supersingular case is identical.  To do this, we check that $\sfn$ is contained in $\sf[\p^n]$ with finite index bounded  independent of $n$.  This suffices to prove the proposition since $\varinjlim \sfn$ is then finite index in $\sf$ and $\sf$ has no proper finite index submodules. (See \cite[Prop 4.14]{Greenberg} for a proof of this fact for elliptic curves over cyclotomic $\Zp$-extensions which generalizes to the case we are considering.  See \cite{Kimfinite} for the supersingular case.  In both cases, one needs as an input the fact that these Selmer groups are $\Lambda$-cotorsion.  This is established in the next section of this paper assuming hypothesis $\ram$ -- see Theorem \ref{thm:bd}.)

Note that 
$$
\sf[\p^n] \subseteq H^1(K_\Sigma/\Ki,A_f)[\p^n] \cong H^1(K_\Sigma/\Ki,A_{f,n})
$$
as $\rhob_f$ is irreducible.  
It a straightforward diagram chase to verify that $\sf[\p^n]$ equals the kernel of
$$
H^1(K_\Sigma/\Ki,A_{f,n}) \to \H_{p,n} \times  \prod_{\substack{\ell \in \Sigma \\ \ell \neq p}} \H_{\ell,n} / \H^{\fin}_{\ell,n} 
$$
where
$$
\H^{\fin}_{\ell,n} = 
\begin{cases} \displaystyle
 A_f^{G_{K_{\ell}}} / \p^n A_f^{G_{K_{\ell}}}
\otimes \Lambda^\vee 
& \text{for~} \ell \text{~inert~in~} K,  \vspace{.2cm}
\\
\displaystyle \prod_{w \in \sigma_{\infty,\ell}} A_f^{G_{K_{\infty,w}}} / \p^n A_f^{G_{K_{\infty,w}}} & \text{for~} \ell \text{~split~in~} K.
\end{cases}
$$
Here, 
we are identifying  $A_f^{G_F} / \p^n A_f^{G_F}$ with its image 
in the exact sequence
$$
0 \to A_f^{G_F} / \p^n A_f^{G_F} \to H^1(F,A_{f,n}) \to H^1(F,A_f)[\p^n] \to 0
$$
for $F=K_\ell$ or $K_{\infty,w}$.

To compare $\sfn$ and $\sf[\p^n]$, we compare the local conditions that define them as subsets of $H^1(K_\Sigma/\Ki,A_{f,n})$.  A  
straightforward computation shows that  for $\ell$ inert in $K$, we have $\H^{\ord}_{\ell,n} = \H^{\fin}_{\ell,n}$. 
Therefore, there is an exact sequence
$$
0 \to \sfn \to \sf[\p^n] \to  \prod_{\substack{\ell \mid N^+ \\ w \in \sigma_{\infty,\ell}}} A_f^{G_{K_{\infty,w}}} / \p^n A_f^{G_{K_{\infty,w}}}.
$$
Another simple computation shows that the final term in this sequence is finite and of size bounded independent of $n$, as desired.
\end{proof}

 \subsection{Evaluating at the trivial character}
The following proposition illustrates how the minimal Selmer group $\sf$ loses information about the rational primes that are inert in $K$ (and thus infinitely split in $\Ki$) while the Greenberg Selmer group $\Sf$ retains information about these primes.
The analogous proposition for cyclotomic extensions is proven in \cite[Section 4]{Greenberg}.  We thus omit a proof here as the arguments of \cite{Greenberg} carry through with only a few changes.

\begin{proposition}
Assume $\sel(K,f)$ is finite.  If $L^{\alg}(\Ki,f)$ (resp.\ $\LL^{\alg}(\Ki,f)$) denotes the characteristic power series of $\sf^\vee$ (resp.\ $\Sf^\vee$) and ${\bf 1}$ denotes the trivial character, then
$$
{\bf 1}\bigl(L^{\alg}(\Ki,f)\bigr) \sim  \left| \sel(K,A_f) \right| \cdot \prod_{v \mid p}  \left|A_f(k_v)\right|^2 \cdot
\prod_{v \mid N^+} \left|k\right|^{t_f(v)}
$$
and
$$
{\bf 1}\bigl(\LL^{\alg}(\Ki,f)\bigr) \sim  \left| \sel(K,A_f) \right| \cdot \prod_{v \mid p}  \left| A_f(k_v) \right|^2 \cdot
\prod_{v \mid N} \left| k \right|^{t_f(v)}.
$$
Here $v$ runs through places of $K$ and  $a \sim b$ if their quotient is a $p$-adic unit.
\end{proposition}



\section{Divisibilities}

\subsection{Statement}
\label{sec:statement}

Recall that $f$ is a normalized newform of weight two and squarefree level $N=N^+N^-$,  $K/\Q$ is a quadratic imaginary field in which all of the prime divisors of $N^+$ (resp.\ $N^-$) are split (resp.\ inert), and $N^-$ is the product of
an odd number of primes.
By \cite[Th\'eor\`eme 2]{Carayol} we may associate to $f$ a 
Galois representation $$\rho_{f} : \GQ \to \GL_{2}(\O_{0})$$
(which gives rise to $T_{f}$ after tensoring with $\O$).
Our goal in this section is to prove the following generalization of the main results of \cite{BD,DI}.  
For a cotorsion $\Lambda$-module $M$, we write
$\chr_{\Lambda}(M)$ for the characteristic ideal of the dual
of $M$.

\begin{theorem} \label{thm:bd}
Assume that $p \geq 5$ and  that
$(\rhob_f,N^-)$ satisfies $\ram$.
If $f$ is $p$-ordinary, then
$\sf$ is a cotorsion $\Lambda$-module and
$$
\chr_\Lambda(\sf) \text{~divides~} \lf  \text{~in~} \Lambda.
$$
If $f$ is $p$-supersingular, $a_{p}=0$, $p$ is split in $K/\Q$  and each prime above $p$ is totally ramified in $\Ki/K$, then
$\sfpm$ is a cotorsion $\Lambda$-module and
$$
\chr_\Lambda(\sfpm) \text{~divides~} \lfpm \text{~in~} \Lambda.
$$
\end{theorem}

\begin{remark}
\label{rmk:hyps}
Bertolini--Darmon \cite{BD} and Darmon--Iovita \cite{DI}
prove this result under two additional hypotheses:
\begin{enumerate}
\item the ring $\O_{0}$ of Fourier coefficients equals $\Zp$,
\item $\eta_f(N^+\!,N^-)$ is a $p$-adic unit.  
\end{enumerate}
We note that our hypothesis  that $(\rhob_f,N^-)$ satisfies $\ram$ is a weakening of this second condition.  Indeed, by level-lowering, if $\eta_f(N^+\!,N^-)$ is a $p$-adic unit, then $\rhob_f$ is ramified at all primes $q \mid N^-$.  (Note that the converse of this statement is not true.)  
Hypothesis \ram~ only demands that if $\rhob_f$ is unramified at some $q \mid N^-$, then $q \not \equiv \pm 1 \pmod{p}$.
\end{remark}

In the remainder of this section we explain how to remove these two hypotheses.
As the necessary changes are identical in the ordinary and supersingular cases, for simplicity we will restrict our attention to the ordinary case.
Our argument is a slight
modification of that of \cite{BD}; we will assume familiarity with
the latter work throughout this section and will focus on the differences.

\subsection{Preparations}

It is essential to our method of proof that we work in a more general setting than the previous section.
To this end, fix $n \geq 1$ and let $f$ denote an ${\O}_{0}/\p_{n}$-valued eigenform for $\Ts$;  for our 
purposes it is most convenient to regard $f$ as a homomorphism
$f : \Ts \to {\O}_{0}/\p_{n}$.  Associated to $f$ and $K$, we have a $p$-adic $L$-function 
$\lf \in \O/\p^n[[\Gamma]]$.  

Let $A_{f}^0$ denote the free $\O_{0}/\p_{n}$-module of rank two endowed with a continuous $\O_{0}/\p_{n}$-linear action of $\GQ$ associated to $f$ by \cite[Th\`eor\'eme 3]{Carayol}.
Define $A_f := A_f^0 \otimes_{\O_0} \O$.
Assuming $\ram$, we may define a residual Selmer group $\sf$ attached to $A_{f}$ and $N^{-}$ as in Section~\ref{sec:modpsel}.
The residual Selmer group $\sf$ is naturally endowed with the structure of a $\Lambda/\p^{n}$-module.  

Finally, we fix an $\O$-algebra 
homomorphism $\varphi : \Lambda \to \O'$ with $\O'$ a 
discrete valuation ring of characteristic zero with maximal ideal $\p'$. 
Let $s_{f}$ denote the
$\O'$-length of $\sf^\vee \otimes_{\Lambda} \O'$ and let
$2t_{f}$ denote the $\O'$-valuation of $\varphi\bp{\lf} 
\in \O'/\varphi(\p)^{n}$ (which we take
to be infinite if $\varphi\bp{\lf}=0$).  
(We use $2t_{f}$ to correspond to the notation of \cite{BD}.)
We will prove the following proposition.

\begin{proposition} \label{prop:bd2}
Fix $n \geq 1$ and $t_{0} \geq 0$.
Let $N^{-}$ denote a squarefree integer divisible by an odd number of primes,
each inert in $K/\Q$.  Let $\tilde{f}$ 
be an ${\O}_{0}/\p_{n+t_{0}}$-valued eigenform 
for $\Ts$ and let $f$ be its projection onto $\O_0/\p_n$.  Assume  that:
\begin{enumerate}
\item The homomorphism $f : \Ts \to {\O}_{0}/\p_{n}$ is
surjective;
\item $(\rhob_f,N^-)$ satisfies hypothesis $\ram$;
\item $2t_{f} < 2t_{0}$.
\end{enumerate}
Then $s_{f} \leq 2t_{f}$. 
\end{proposition}

We claim that this proposition implies Theorem~\ref{thm:bd}.  Indeed,
fix $f$ as in Theorem~\ref{thm:bd} and a homomorphism $\varphi$ as
above.  If $\varphi\bp{\lf}=0$, then certainly
$\varphi\bp{\lf}$ belongs to the $\O'$-Fitting ideal of
$\sf^{\vee} \otimes_{\Lambda} \O'$.  Otherwise, choosing $t_{0}$
larger than the $\O'$-valuation of $\varphi\bp{\lf}$ and applying
Proposition~\ref{prop:bd2} for all $n$, we obtain again that
$\varphi\bp{\lf}$ belongs to the Fitting ideal of
$\sf^{\vee} \otimes_{\Lambda} \O'$.  
Applying this for all $\varphi$, Proposition~\ref{prop:limit} and
\cite[Proposition 3.1]{BD} yield Theorem~\ref{thm:bd}.

\subsection{Construction of cohomology classes}

Fix $N^{-}$ divisible only by primes inert in $K/\Q$ and a surjective
homomorphism $f : \Ts \to \O_{0}/\p_{n}$ such that
the residual representation $\rhob_{f}$ satisfies hypothesis $\ram$.  Write $T^0_{f}$ for the Galois
representation associated to $f$ over $\O_{0}/\p_{n}$ and set
$T_{f} = T^0_{f} \otimes_{\O_{0}} \O$.

We say that a rational prime $\ell$ is {\it admissible} relative to $f$ if:
\begin{enumerate}
\item $\ell$ does not divide $N^{-}$;
\item $\ell$ is inert in $K/\Q$;
\item $p$ does not divide $\ell^{2}-1$;
\item One of $\ell + 1 - f(T_{\ell})$ or $\ell + 1 + f(T_{\ell})$ equals zero
in $\O_{0}/\p_{n}$.
\end{enumerate}

Fix an admissible prime $\ell$ relative to $f$.
In this section we give the construction of a cohomology class
$$\kappa(\ell) \in \hat{H}^{1}(\Ki,T_{f}) = \varprojlim_n H^1(\Ki,T_f/\p^nT_f)$$ 
which is central to the argument.

For any admissible prime $\ell$, the arguments of
\cite[Theorem 5.15]{BD}
yield a surjective homomorphism
$$f_{\ell} : \Tv{N^+\!}{N^-\ell} \to{\O}_{0}/\p_{n}$$
which agrees with $f$ at all Hecke operators away from $\ell$ and which sends the Atkin-Lehner operator
$U_{\ell}$ to the unique $\epsilon \in \{ \pm 1 \}$ such that
$\ell + 1 - \epsilon f(T_{\ell})$ vanishes in $\O_{0}/\p_{n}$. 
(This construction uses the freeness of a certain character group which is  established under hypothesis $\ram$ in Theorem \ref{thm:free} later in the paper.)
Note that the
surjectivity of $f_{\ell}$ is automatic from the surjectivity of $f$ and the admissibility of $\ell$.  

Write $\Ifl$ for the kernel of $f_{\ell}$,
$\Jl$ for the Jacobian of the Shimura curve $X_0(N^+\!,N^{-}\ell)$, and
$\Ta_{p}(\Jl)$ for its $p$-adic Tate module.
The key result, which corresponds to \cite[Theorem 5.17]{BD}, is the following.

\begin{proposition}
\label{prop:gal}
The Galois representations $\Ta_{p}(\Jl)/\Ifl$ and $T_{f}^0$ are
isomorphic.
\end{proposition}
\begin{proof}
As the Frobenius traces of $\GQ$ on these two representations agree for
all primes away from $N\ell$, it suffices to show that
$T:=\Ta_{p}(\Jl)/\Ifl$ is free of rank two over ${\O}_{0}/\p_{n}$.
Since $\p_{1}$ is the maximal ideal of ${\O}_{0}$,
the first step of \cite[Theorem 5.17]{BD} applies to show that
$T/\p_{1}$ is two-dimensional over $k_{0}$.  We will use this to deduce the desired result for $T$.

We first show that $T$ has a free ${\O}_{0}/\p_{n}$-submodule of rank one.
Let $\Phi_{\ell}$ denote the group of connected components of the
N\'eron model of $\Jl$ over the Witt vectors of $\F_{\ell^{2}}$.
The proof of \cite[Lemma 5.15]{BD} shows that $\Phi_{\ell}/\Ifl$ is isomorphic to ${\O}_{0}/\p_{n}$.  Let
$c \in \Phi_{\ell}/\Ifl$ correspond to $1 \in{\O}_{0}/\p_{n}$ under some
such isomorphism.  Applying the argument of
\cite[Lemma 5.16]{BD} to  $c$ yields an integer $n'$ and an element
$t \in \Jl[p^{n'}](\Qlur)/\Ifl$ which maps onto $c$
under the natural map
$$\Jl[p^{n'}](\Qlur)/\Ifl \to \Phi_{\ell}/\Ifl.$$
Since this map respects the Hecke actions and thus 
is ${\O}_{0}/\p_{n}$-linear, the cyclic ${\O}_{0}/\p_{n}$-module generated
by $t$ surjects onto a free ${\O}_{0}/\p_{n}$-module of rank one and thus must itself
be free of rank one.  As
$$\Jl[p^{n'}](\Qlur)/\Ifl \hookrightarrow T$$
this yields the desired submodule.

As $\GQ$
acts irreducibly on $T/\p_{1}$, we may choose an element
$g \in \GQ$ so that $t$ and $u := gt$ are a basis of
$T/\p_{1}$; here $t$ is the element constructed above generating a
free ${\O}_{0}/\p_{n}$-module of rank one.
Note that $u$ also generates a free ${\O}_{0}/\p_{n}$-module 
of rank one.

We will show that $t$ and $u$ are an ${\O}_{0}/\p_{n}$-basis of
$T$.  They span by Nakayama's lemma, so it suffices to show
that in any relation
\begin{equation} \label{eq:o1}
\alpha t = \beta u
\end{equation}
with $\alpha,\beta \in{\O}_{0}$, we must have $\alpha,\beta \in \p_{n}$.
Note that we must have $\alpha,\beta \in \p_{1}$ since $t,u$ are a basis
modulo $\p_{1}$.  In fact, if $d$
denotes the minimal $\O$-valuation of any element of $\p_{1}$, then
we have that $\p_{1} = \p_{d}$, so that
\begin{equation} \label{eq:o2}
\alpha,\beta \in \p_{d}.
\end{equation}
Since $\GQ$ surjects onto $\GL_{2}(k_{0})$, we may choose an element
$h \in \GQ$ with the property that
$$ht = 2t + t'; \qquad hu = u + u'$$
where
\begin{equation} \label{eq:o3}
t',u' \in \p_{1}T = \p_{d}T.
\end{equation}
Applying $h-1$ to (\ref{eq:o1}), we find that
\begin{equation} \label{eq:o4}
\alpha t =  \beta u' - \alpha t'.
\end{equation}

Suppose now that we know that $\alpha,\beta \in \p_{r}$ for some $r < n$.
We will show that in fact $\alpha,\beta \in \p_{r'}$ for some $r' > r$;
the fact that $\alpha,\beta \in \p_{n}$ then 
follows by induction.  By (\ref{eq:o3}) and (\ref{eq:o4}) we have
that
\begin{equation} \label{eq:o5}
\alpha t \in \p_{r+d}T.
\end{equation}
If $r+d \geq n$, then $\alpha t = 0$; since $t$ generates
a free ${\O}_{0}/\p_{n}$-module, it follows that $\alpha \in \p_{n}$.  The
same argument shows that $\beta \in \p_{n}$ as well, so that in this case
we are done.

If $r+d < n$, then
multiplying both sides of (\ref{eq:o5}) by $\p_{n-r-d}$ and using that
$\p_{n-r-d}\p_{r+d} \subseteq \p_{n}$ yields
$$\alpha \p_{n-r-d}t = 0.$$
Since $t$ generates a free $\O_{0}/\p_{n}$-module, it follows that
$$\alpha \p_{n-r-d} \subseteq \p_{n}.$$
Let $a,b > 0$ be such that
$$\p_{n-r-d}\O = \p^{a}; \qquad \p_{n}\O = \p^{b}.$$
Then we have that $\alpha \in \p^{b-a} \cap \O_0 = \p_{b-a}$.
Clearly $b \geq n$.  Also, $a < n-r$ since $\O_0$ contains elements of
all valuations which are multiples of $d$ and some multiple of $d$
lies between $n-r-d$ and $n-r$.  It follows that $b-a > r$, so that
$\alpha \in \p_{r'}$ for $r' = b -a > r$.  Since by (\ref{eq:o1}) we have
$\beta u \in \p_{r+d}T$ as well, an identical argument shows that
$\beta \in \p_{r'}$, as desired.  This completes the proof.
\end{proof}

With this result in hand, the construction of the cohomology class
$$\kappa_{0}(\ell) \in \hat{H}^{1}(\Ki,T_{f}^0)$$
proceeds as in \cite[Sections 6,7]{BD}.  
Defining $\kappa(\ell)$ as the image of $\kappa_{0}(\ell)$ under the
natural map
$$
\hat{H}^{1}(\Ki,T_{f}^0) \to \hat{H}^{1}(\Ki,T_{f})
$$
the proof of the two explicit reciprocity laws
\begin{align*}
\delta_{\ell}(\kappa(\ell)) &= \lf  ~\text{~~~in~}\Lambda/\p_n \\
v_{\ell_{2}}(\kappa(\ell_{1})) &= \lgn ~\text{~~~in~}\Lambda/\p_n
\end{align*}
proceeds as in \cite[Sections 8--9]{BD}.  

\subsection{Euler system arguments}

We now give the proof of Proposition~\ref{prop:bd2} via a modification of
the Euler system arguments of \cite[Section 4]{BD}.
Our proof proceeds by induction on $t_{f}$.  Let $f$ be
an eigenform as in Proposition~\ref{prop:bd2}.
When $t_{f}=0$, the proof of \cite[Proposition 4.7]{BD} carries over to
prove that $\sf$ is trivial, as required.  Assume therefore that
$t_{f} > 0$.

As in \cite{BD}, for any $(n+t_{f})$-admissible prime $\ell$ we may
construct from $\kappa(\ell)$ and an $(n+t_{f})$-admissible set of primes 
$S$ a cohomology class 
$$\kappa'_{\varphi}(\ell) \in \hat{H}^{1}_{S}(\Ki,T_{f})$$
satisfying \cite[Lemma 4.5 and 4.6]{BD}.  
Let
$\Pi$ denote the set of $(n+t_{0})$-admissible primes $\ell$ for which
$\ord_{\p'}(\kappa_{\varphi}(\ell))$ is minimal.  The set $\Pi$ is
non-empty by \cite[Theorem 3.2]{BD}, and, writing $t$ for the value
$\ord_{\p'}(\kappa_{\varphi}(\ell))$ for $\ell \in \Pi$, by
\cite[Lemma 4.8]{BD} we have $t < t_{f}$.
(Note that we are using $(n+t_{0})$-admissible primes, rather than
$(n+t_{f})$-admissible primes as in \cite{BD}.  This is necessary to
facilitate our induction but has no effect on the results used above.)

Fix $\ell_{1} \in \Pi$ and let $s \in H^{1}(K,T_{f}) \otimes \O'/\p'$ 
denote the image of
$\kappa'_{\varphi}(\ell)$ in
$$\hat{H}^{1}_{S}(\Ki,T_{f}) \otimes \O'/\p' \subseteq \hat{H}^{1}_{S}(\Ki,T_{f})/\m_{\Lambda} \otimes \O'/\p'
\subseteq H^{1}(K,T_{f}) \otimes \O'/\p'.$$
By \cite[Theorem 3.2]{BD} there exists an
$(n+t_{0})$-admissible prime $\ell_{2}$ such that
$v_{\ell_{2}}(s) \neq 0$; here
$$v_{\ell_{2}} : \hat{H}^{1}(K_{\infty,\ell_{2}},T_{f}) \to
\hat{H}^{1}_{\text{fin}}(K_{\infty,\ell_{2}},T_{f})$$
is as in \cite{BD}.
Note that
\begin{equation} \label{eq:bd1}
t = \ord_{\p'}(\kappa_{\varphi}(\ell_{1})) \leq
\ord_{\p'}(\kappa_{\varphi}(\ell_{2})) \leq
\ord_{\p'}(v_{\ell_{1}}(\kappa_{\varphi}(\ell_{2}))).
\end{equation}
(Here the first inequality is by the definition of $\Pi$ and the
second follows from the fact that $v_{\ell_{1}}$ is a homomorphism.)
However, by \cite[Corollary 4.3]{BD} we have
$$\ord_{\p'}(v_{\ell_{1}}(\kappa_{\varphi}(\ell_{2}))) = \ord_{\p'}(v_{\ell_{2}}(\kappa_{\varphi}(\ell_{1}))).$$
Furthermore, since $v_{\ell_{2}}(s) \neq 0$, we must have
$$\ord_{\p'}(v_{\ell_{2}}(\kappa_{\varphi}(\ell_{1}))) = \ord_{\p'}(\kappa_{\varphi}(\ell_{1})).$$
It follows that the inequalities in (\ref{eq:bd1}) must be equalities;
in particular,
$$\ord_{\p'}(\kappa_{\varphi}(\ell_{2})) = t,$$
so that $\ell_{2} \in \Pi$.

Let $g$ denote the $\O_{0}/\p_{n+t_{0}}$-valued eigenform for
$\Tv{N^+\!}{N^{-}\ell_{1}\ell_{2}}$ attached to $f$ and $(\ell_{1},\ell_{2})$
by \cite[Proposition 3.12]{BD}.  By \cite[Theorem 4.2]{BD} we have
$$v_{\ell_{2}}(\kappa(\ell_{1})) = \lgn.$$
Thus $t_{g} = t < t_{f}$.  The eigenform $g$ satisfies all of the
hypotheses of Proposition~\ref{prop:bd2}, so that we may now apply
the induction hypothesis to conclude that $s_{g} \leq 2t_{g}$.
From here one argues as in \cite[pp.\ 34--35]{BD} to conclude that
$s_{f} \leq 2t_{f}$ as well.  This completes the proof.

\section{Algebraic $\mu$-invariants}

We return now to the notation of Section 3.
For simplicity we focus initially on the ordinary case and state the
results in the supersingular case at the end of the section; the proofs
are identical (using \cite[Proposition 4.16]{IP} to check hypothesis
(\ref{4}) of Proposition~\ref{prop:surj}).

Our comparison of $\mu$-invariants of Selmer groups relies crucially on the exactness of the following sequences.

\begin{proposition} \label{prop:lg}
Assume that $p \geq 5$ and that $(\rhob_f,N^-)$ satisfies $\ram$.
The defining sequences
$$0 \to \sf \to H^{1}(K_{\Sigma}/\Ki,A_{f}) \to
\prod_{\ell \in \Sigma} \H_{\ell} \to 0$$
$$0 \to \Sf \to H^{1}(K_{\Sigma}/\Ki,A_{f}) \to
\H_p \times \prod_{\ell \in \Sigma^+} \H_\ell \times \prod_{\ell \in \Sigma^-}
\H_{\ell}/\H_{\ell}^{\un} \to 0$$
are exact.
\end{proposition}

\begin{proof}
In Appendix A, we include a general proposition on the surjectivity of global-to-local maps. To check the hypotheses of this proposition, note that the first is immediate, the second follows from Theorem \ref{thm:bd},  the third is a consequence of the irreducibility of $V_{f}$, and the last follows from the fact that $\sum_{\p \mid p} r_{\p} = 2$ (see \cite[Proposition 1]{Gord}).
\end{proof}

\begin{corollary}
\label{cor:mucomp}
Assume $(\rhob_f,N^-)$ satisfies $\ram$.  Then the sequence
$$0 \to \sf \to \Sf \to \prod_{\ell \mid N^-} \H_{\ell}^{\un} \to 0$$
is an exact sequence of cotorsion $\Lambda$-modules.  In particular,
$$\lambda\bp{\sf} = \lambda\bp{\Sf}$$
and
$$\mb{\Sf} = \mb{\sf} + \sum_{\ell \mid N^-}
t_f(\ell).$$
\end{corollary}
\begin{proof}
This is immediate from Proposition~\ref{prop:lg} and Lemma~\ref{lemma:tam}.
\end{proof}

Combining this corollary with the results of \cite{Vmu}, we thus obtain the following
theorem, which we state in both the ordinary and supersingular cases.

\begin{theorem}
\label{thm:mualg}
Let $f$ be a normalized newform of weight two and squarefree level $N=N^+N^-$ with $N^-$ divisible by an odd number of primes.  Assume hypothesis $\ram$ holds for $(\rhob_f,N^-)$.
\begin{enumerate}
\item If $f$ is $p$-ordinary, then
$$\mb{\sf}=0 \text{~and~}
\mb{\Sf}=\sum_{\ell \mid N^{-}} t_f(\ell).$$
\item If $f$ is $p$-supersingular, $a_{p}=0$, $p$ is split in $K/\Q$ and each prime above $p$ is totally ramified in $\Ki/K$, then
$$\mb{\sfpm}=0 \text{~and~} \mb{\Sfpm}=\sum_{\ell \mid N^{-}} t_f(\ell).$$
\end{enumerate}
\end{theorem}

\begin{proof}
By Theorem \ref{thm:muanal}.1, we know that 
$\mb{\lf}$ vanishes.  Thus the $\O[[\Gamma]]$-divisibility of 
Theorem \ref{thm:bd}  implies that $\mb{\sf}$ vanishes as well.
Corollary \ref{cor:mucomp} then gives the value of $\mb{\Sf}$.  The supersingular case is identical.
\end{proof}

\section{The $\mu$-part of the main conjecture}
\label{sec:mumc}

\subsection{The main conjecture}

The anticyclotomic Iwasawa main conjecture for a modular form $f$ relates its $p$-adic $L$-function to the characteristic ideal of its Selmer group.  As we have seen, there are two choices of each object in this setting.  The main conjecture predicts that they correspond as follows.  

\begin{conjecture}  \label{conj:mc}
Let $f$ be a modular form of weight two and squarefree level $N=N^+N^-$ with $N^-$ the product of an odd number of primes.  Assume that the residual representation $\rhob_{f}$ is absolutely irreducible. 
\begin{enumerate}
\item If $f$ is $p$-ordinary, then $\sf$ and $\Sf$ are $\Lambda$-cotorsion,
$$\chr_\Lambda\bp{\sf} = \lf \cdot \Lambda,$$
$$\chr_\Lambda\bp{\Sf} = \Lf \cdot \Lambda.$$
\item If $f$ is $p$-supersingular, $\sfpm$ and $\Sfpm$ are $\Lambda$-cotorsion,
$$\chr_\Lambda\bp{\sfpm} = \lfpm \cdot \Lambda,$$
$$\chr_\Lambda\bp{\Sfpm} = \Lfpm \cdot \Lambda.$$
\end{enumerate}
\end{conjecture}

The formulae on $\mu$-invariants in Theorems~\ref{thm:muanal}, \ref{thm:muanalss} and~\ref{thm:mualg} immediately yield the $\mu$-part of this conjecture for $\sf$ and $\lf$.  However, the corresponding formulae for the $\mu$-invariants of $\Sf$ and $\Lf$ do not immediately appear identical. Indeed, 
by Theorems~\ref{thm:muanal} and~\ref{thm:mualg},
this equality of $\mu$-invariants reduces to the equality
\begin{eqnarray}
\label{eqn:mumu}
\ord_\p \left(\frac{\cf}{\xf }\right)  =
\sum_{\ell \mid N^-} t_f(\ell).
\end{eqnarray}
The right hand side of (\ref{eqn:mumu}) is a purely local expression while the left hand side is a difference of global terms.  

In the remainder of this section,
we will reduce (\ref{eqn:mumu}) to an equality involving degrees of modular parameterizations arising from Shimura curves which was established by Ribet and Takahashi \cite{RT,T}.

\subsection{Character groups}
\label{sec:char}

For this subsection, we fix a factorization $N=N_1 N_2$ such that $N_2$ has an {\it even} number of prime divisors.  (In the following subsection, we will take $N_1 = N^+r$ and $N_2 = N^-/r$ for some $r \mid N^-$.)  Let $J = J_0(N_1,N_2)$ denote the Jacobian of the Shimura curve of level $N_1$ attached to the indefinite quaternion algebra ramified exactly at the primes dividing $N_2$.  

The special fiber at $r \mid N$ of the N\'eron model of $J$ is an extension of an abelian variety by a torus.  Let $\XJr = \X_r(N_1,N_2)$ denote the $\O$-completion of the character group of this torus endowed with its natural action of $\T_0(N_1,N_2)$.
If $r \mid N_1$, then this action factors through the $r$-new quotient of this Hecke algebra and, moreover,  $\XJr$ is a faithful $\T_0(N_1/r,rN_2)$-module (see, for instance, \cite[Proposition 5.8.3]{BD}).
Let $\m_f \subseteq \T_0(N_1,N_2)$ be the maximal ideal corresponding to $f$; set $\hT_0(N_1,N_2) = \T_0(N_1,N_2)_{\m_f} \otimes_{\Zp} \O$ and $\hXJr = \XJr_{\m_f}$.

\begin{theorem}
\label{thm:free}
If $(\rhob_f,rN_2)$ satisfies hypothesis $\ram$, then
$\hXJr$ is free of rank one over the Hecke algebra $\hT_0(N_1/r,rN_2)$.
\end{theorem}

\begin{proof}
By \cite[Theorem 2.1]{Wiles}, $J_0(N)[\m_f]$ has dimension 2 over $k_0$.  Since $(\rhob_f,N_2)$ satisfies hypothesis $\ram$, by \cite[Corollary 8.11, Remark 8.12]{Helm}, $J_0(N_1,N_2)[\m_f]$ has dimension 2 over $k_0$.   Since $(\rhob_f,r)$ satisfies hypothesis $\ram$, a standard application of Mazur's principal (for instance, \cite[Lemma 6.5]{Helm}) implies $\XJr / \m_f \XJr$ has dimension 1 over $k$. Thus, by Nakayama's lemma and the faithfulness of the Hecke-action, we deduce that $\hXJr$ is free of rank one over $\hT_0(N_1/r,rN_2)$.
\end{proof}

Consider now the optimal quotient $\xi : J \to A$ attached to $f$; thus $A$ is an abelian variety and $\ker(\xi)$ is connected.
Let $\hXAr$ and $\hat{\X}_r(A^\vee)$ be the analogues of $\hXJr$ for  $A$ and its dual abelian variety $A^\vee$.
As in \cite[p.\ 208]{Khare}, we may fix an isomorphism $\hXAr \cong \hat{\X}_r(A^\vee)$.  With this isomorphism fixed, the map $\xi$ induces maps 
$$
\xi^*: \hXJr \to \hXAr \text{~~and~~} \xi_*: \hXAr \to \hXJr.
$$
Moreover, as in \cite[p.\ 207]{Khare}, $\xi^* \xi_*$ acts on $\hXAr$ by multiplication by some element $\delta_f(N_1,N_2)$ in $\O$.  Changing our chosen isomorphism above will only change this number by a $p$-adic unit and so the $\O$-ideal $\bp{\delta_f(N_1,N_2)}$ is well-defined.  We simply write $\delta_f(N)$ for $\delta_f(N,1)$.

Let $\PJr$ denote the component group of $J$ at $r$ and set $\hPJr = (\PJr \otimes_\Z \O)_{\m_f}$.  We define $\hPAr$ analogously.
We state here two propositions summarizing the properties of these character and component groups that will be needed in what follows.

\begin{proposition}\label{prop:comp}
~
\begin{enumerate}
\item The monodromy pairings $\langle \cdot, \cdot \rangle_{A}$
and $\langle \cdot , \cdot \rangle_J$
induce exact sequences
$$
0 \to \hXAr \to \hXAr^\vee \to \hPAr \to 0;
$$
$$
0 \to \hXJr \to \hXJr^\vee \to \hPJr \to 0.
$$
\item If $r \mid N_2$, then $\hPJr = 0$ and $\hPAr$ is $\O$-cyclic of order  $\left|k\right|^{t_f(r)}$.  
\end{enumerate}
\end{proposition}

\begin{proof}
The first part is \cite[Theorem 11.5]{Groth}; the second part follows from \cite[Proposition 3]{Khare}.
\end{proof}

Let $\hXJr^f \subseteq \hXJr$ denote the subgroup on which 
$\hT_{0}(N_{1},N_{2})$
acts via $\pi_f$; in particular
$\xi^* \hXAr \subseteq \hXJr^f$. Then $\hXJr^f$ is a free $\O$-module of rank 1; let $g_r$ denote a generator of this module.

\begin{proposition}
\label{prop:char}Let $r \mid N_1$.
\begin{enumerate}
\item 
\label{prop:char.1}
$\hXJr^f / \xi^* \hXAr$ is $\O$-cyclic with size $|k|^{t_f(r)}$.
\item If $\hXJr$ is free over $\hT_0(N_1/r,rN_2)$, then
$$
\bp{\langle g_r, g_r \rangle_J} =  \bp{\eta_f(N_1/r,rN_2)}
$$
\end{enumerate}
\end{proposition}

\begin{proof}
The first part follows from \cite[Lemma 2]{Khare}.  For the second part, by Proposition \ref{prop:comp}, the monodromy pairing on $\hXJr$ is perfect.  Thus, since $\hXJr$ is free over $\hT_0(N_1/r,rN_2)$,
by \cite[Lemma 4.17]{DDT}, $\langle g_r, g_r \rangle_J$ computes the congruence number $\eta_f(N_1/r,rN_2)$.
\end{proof}

\subsection{Comparing definite and indefinite Shimura curves}

Fix a divisor $r$ of $N^-$
and let $\XJr = \X_r(N_1,N_2)$ be the character group of the previous subsection where $N_1 = N^+r$ and $N_2 = N^-/r$.  Let $X_{N^+\!,N^-}$ denote the Shimura curve (as in Section \ref{sec:analytic}) of level $N^+$ attached to the definite quaternion algebra ramified at the primes dividing $N^-$.

\begin{proposition}
\label{prop:compare}
For each $r \mid N^-$, there is a canonical Hecke-equivariant isomorphism
$$
\Pic(X_{N^+\!,N^-}) \otimes \O \cong \X_r(N^+r,N^-/r).
$$
Moreover, this isomorphism takes the intersection pairing on $\Pic(X_{N^+\!,N^-})$ to the monodromy pairing on $\X_r(N^+r,N^-/r)$.
\end{proposition}

\begin{proof}
See \cite[Theorem 4.3]{Kohel}.
\end{proof}

We now give the proof of Lemma \ref{lemma:unit} of Section \ref{sec:analytic}.

\begin{proof}[Proof of Lemma \ref{lemma:unit}]

Let $r$ be any divisor of $N^-$.  By Proposition \ref{prop:compare}, we have 
$$
\M = \Pic(X_{N^+\!,N^-}) \otimes \O \cong \X_r(N^+r,N^-/r).
$$ 
Thus, the lemma follows from \cite[Lemma 2.2, Theorem 2.7]{T} which proves the analogous statement for the character group $\X_r(N^+r,N^-/r)$.
\end{proof}

\subsection{Modular degrees}

Throughout this section, for $x \in \O$, we write $\bp{x}$ for the $\O$-ideal generated by $x$.  The following proposition (and its proof) is essentially \cite[Theorem 2.3]{T}.

\begin{proposition}
\label{prop:pairdelta}
We have 
$$
\ord_\p \bp{\delta_f(N_1,N_2)} = t_f(r) +  \ord_\p \bp{\langle g_r, g_r \rangle_J}
$$
for $r \mid N_1$.
\end{proposition}

\begin{proof}
Let $g_r$ be a generator of $\hXJr^f$.
By Proposition \ref{prop:char}, if $x_r$ is a generator of $\hXAr$, then $\xi_*(x_r) = c_f(r) g_r$ for some $c_f(r) \in \O$ such that $\ord_\p c_f(r) = t_f(r)$.  
Thus,
\begin{align*}
\delta_f(N_1,N_2) \cdot \langle x_r , x_r \rangle_A 
&= \langle x_r , \xi^* \xi_* x_r \rangle_A 
= \langle \xi_* x_r ,  \xi_* x_r \rangle_J \\
&= \langle c_f(r) g_r  ,  c_f(r) g_r \rangle_J 
= c_f(r)^2 \cdot \langle g_r  , g_r \rangle_J.
\end{align*}
By Proposition \ref{prop:comp},
$\ord_\p \bp{\langle x_r , x_r \rangle_A} = t_f(r)$ which proves the proposition. 
\end{proof}

The following proposition establishes the equality of the $p$-parts of a congruence number and a degree of a modular parameterization in the case of a modular curve, that is, $N_2=1$.  For more general results along these lines see \cite{AU,CK,ARS}.

\begin{proposition}
\label{prop:deltaeta}
We have, as $\O$-ideals,
$$
\bp{\delta_f(N)} = \bp{\eta_f(N)}.
$$
\end{proposition}

\begin{proof}
By level-lowering, there exists some prime $r \mid N$ such that $\rhob_f$ is ramified at $r$. Thus, the Tamagawa exponent $t_f(r)$ is zero and by
Proposition \ref{prop:pairdelta} we have
$$
\bp{\delta_f(N)} =   \bp{\langle g_r, g_r \rangle_J}
$$
as ideals of $\O$ where $g_r$ is a generator of $\hat{\X}_r(N/r,r)^f$.

Since $\rhob_f$ is ramified at $r$, it is immediate that $(\rhob_f,r)$ satisfies \ram~ and, thus, by Theorem \ref{thm:free}, $\hXJr$ is free over $\hT_0(N/r,r)$.  By Proposition \ref{prop:char} it follows that
$$
 \bp{\langle g_r, g_r \rangle_J} = \bp{ \eta_f(N/r,r) }.
$$
But, as $\rhob_f$ is ramified at $r$,
$$
\bp{\eta_f(N)} = \bp{ \eta_f(N/r,r) }
$$
as any form congruent to $f$ is automatically new at $r$.  Combining these equalities yields the proposition.
\end{proof}

\subsection{The $\mu$-part of the main conjecture}

\begin{theorem}
\label{thm:formula}
Let $f$ be a modular form of weight two and squarefree level $N=N^+ N^-$
with $N^-$ the product of an odd number of primes.  Assume that
$(\rhob_{f},N^-)$ satisfies hypotheses $\ram$.  Then
$$
\ord_\p \left(\frac{\cf}{\xf }\right)  =
\sum_{\ell \mid N^-} t_f(\ell).
$$  
\end{theorem}

\begin{proof}
Recall that $\xf = \langle g_f, g_f \rangle$ where $g_f$ is a generator of $\M^f$, the subspace of $\Pic(X_{N^+\!,N^-}) \otimes_{\Z} \O$ where $\T_{0}(N^+\!,N^-)$ acts via $\pi_f$.  By Proposition \ref{prop:compare}, for any $r \mid N^-$, we have $\langle g_f , g_f \rangle = \langle g_r , g_r \rangle_J$ where $g_r$ is a generator of $\hXJr^f$.  Thus, Proposition \ref{prop:pairdelta} and  Proposition \ref{prop:deltaeta} yield
$$
\ord_\p \left(\frac{\cf}{\xf}\right)  =
\ord_\p \left(\frac{\delta_f(N,1)}{\delta_f(N^+r,N^-/r)} \right) + t_f(r).
$$
The main result of Ribet--Takahashi \cite{RT,T} is
$$
\ord_\p \left(\frac{\delta_f(N,1)}{\delta_f(N_1,N_2)} \right) 
= \sum_{r \mid N_2} t_f(r)
$$
for any factorization $N=N_1 N_2$ where $N_2$ has an even number of prime factors.  
Combining these two equalities then yields the theorem.
\end{proof}

\begin{theorem}
\label{thm:mu}
Assume $(\rhob_f,N^-)$ satisfies $\ram$. 
\begin{enumerate}
\item
If $f$ is $p$-ordinary, then 
$$\mb{\sf}=\mb{\lf}=0 \text{~and}$$
$$\mb{\Sf}=\mb{\Lf}.$$
\item 
If $f$ is $p$-supersingular, $a_p=0$, $p$ is split in $K/\Q$  and each prime above $p$ is totally ramified in $\Ki/K$, then
$$\mb{\sfpm} = \mb{\lfpm} = 0 \text{~and}$$
$$\mb{\Sfpm} = \mb{\Lfpm}.$$
\end{enumerate}
\end{theorem}

\begin{proof}
The equality for the minimal Selmer group and $\lf$ follows from Theorem \ref{thm:muanal} and Theorem \ref{thm:mualg} as all of these $\mu$-invariants are zero.  The equality for the Greenberg Selmer group follows from these theorems and from Theorem \ref{thm:formula} as 
$$
\mb{\Lf} = \ord_\p \left( \frac{\eta_f(N)}{\xf}\right) = \sum_{r \mid N^-} t_f(r) = \mb{\Sf}.
$$
The supersingular case follows identically by appealing to Theorem \ref{thm:muanalss}.
\end{proof}

\subsection{Quantitative level-lowering}
\label{sec:pf}

When $(\rhob_f,N^-)$ satisfies \ram, we have
$$
\xf = \eta_f(N^+\!,N^-) 
$$ 
by Theorem \ref{thm:free} and Proposition \ref{prop:char}.  
In this case, the $\mu$-part of the main conjecture is thus equivalent to the equality
\begin{eqnarray}
\label{eqn:formula}
\ord_\p \left( \frac{\eta_f(N)}{\eta_f(N^+\!,N^-)}\right) = \sum_{r \mid N^-} t_f(r).
\end{eqnarray}
We propose the following formula, which would explain this equality in general.

\begin{quote}
Let $f$ be a weight two eigenform of squarefree level $N=aqb$ with $q$ prime.  Then
$$
\ord_\p \bp{\cng{f}{aq}{b}} = t_f(q) +  \ord_\p \bp{\cng{f}{a}{q b}}.
$$
as ideals of $\O$.
\end{quote}
\vspace{.1cm}

\begin{remark}~
\begin{enumerate}
\item Equation (\ref{eqn:formula}) follows immediately from this formula as
\begin{align*}
\ord_\p \bp{\cf} &= t_f(q) +  \ord_\p \bp{\cng{f}{N/q}{q}} \\
&= t_f(q) + t_f(q') + \ord_\p \bp{\cng{f}{N/(q q')}{q q'}} \\
&= \dots = \sum_{q \mid N^-} t_f(q) + \ord_\p \bp{\cng{f}{N^+\!}{N^-}}.
\end{align*}

\item The formula of Ribet and Takahashi \cite{RT} involving degrees of modular parameterization arising from Shimura curves is completely analogous to the above formula when one changes the new-part of the level two primes at a time.

\item The analogous formula for congruence numbers arising from changing the new-part of the level by two primes was used by Khare \cite{Khare} to relate certain new quotients of a Hecke algebra to the full Hecke algebra in a case where hypothesis \ram~was satisfied.
\end{enumerate}
\end{remark}

This formula can be regarded as a quantitative form of level-lowering much like Wiles' numerical criterion \cite{Wiles} can be viewed as a quantitative version of level-raising.
More precisely, let $f$ be an eigenform of weight two and level $N$, and suppose that $\rhob_{f}$ is unramified at a prime $q \mid N$.  By
definition, $t_{f}(q) > 0$, so that the proposed formula implies that $\eta_{f}(N,1)/\eta_{f}(N/q,q)$ is a non-unit.  That is, there exists an eigenform $g$ of level $N$ which is congruent to $f$, but which is old at $q$.  This is precisely what is  predicted by level-lowering.

We conclude this section with a proof of this formula assuming $\ram$. Note that  this hypothesis puts us into the case where ``Mazur's principle" applies to establish level-lowering.

\begin{theorem}
\label{thm:conj}
Let $f$ be a newform of weight two and squarefree level $N=aqb$.  Assume $(\rhob_f,bq)$ satisfies $\ram$ and that there are at least two primes at which $\rhob_f$ is ramified.  Then 
$$
\ord_\p \bp{\cng{f}{aq}{b}} = t_f(q) +  \ord_\p \bp{\cng{f}{a}{q b}}.
$$
\end{theorem}

\begin{proof}
Note that if $\rhobar_f$ is ramified at $r \mid N$, then
$$
\bp{\cng{f}{cr}{d}} = \bp{\cng{f}{c}{rd}}
$$
for $N=crd$ since any form congruent to $f$ must be $r$-new.
Since $\rhob_f$ is ramified at two distinct primes, using the above observation, we may assume that $b$ has an even number of prime factors and that there is some prime $r \mid a$ at which $\rhob_f$ is ramified.

On the one hand, applying Proposition \ref{prop:pairdelta} at the prime $r$ (which is valid since $b$ has an even number of prime factors) yields
\begin{align*}
\ord_\p \bp{\delta_f(aq,b)} &= 
t_f(r) + \ord_\p \bp{\langle g_r, g_r \rangle_J} \\
&= t_f(r) + \ord_\p \bp{\eta_f(aq/r,rb)} 
= \ord_\p \bp{\eta_f(aq,b)}.
\end{align*}
The second equality follows from Proposition \ref{prop:char} and Theorem \ref{thm:free} (as we are assuming \ram).  The third equality follows 
since $\rhobar_f$ is ramified at $r$.

On the other hand, applying Proposition \ref{prop:pairdelta} at the prime $q$ yields
$$
\ord_\p \bp{\delta_f(aq,b)} = t_f(q) + \ord_\p \bp{\eta_f(a,q b)}. 
$$
Hence
$$
\ord_\p \bp{\eta_f(aq,b)} = t_f(q) + \ord_\p  \bp{\eta_f(a,q b)}
$$
proving the theorem.
\end{proof}

\section{$\lambda$-invariants and congruences}

In the papers \cite{GV,EPW} explicit formulae are given for the differences of $\lambda$-invariants of cyclotomic Selmer groups of congruent modular forms.  These results transfer verbatim to the setting of anticyclotomic Selmer groups for modular forms $f$ such that $N(\rho_f)/N(\rhob_f)$ is only divisible by primes split in $K$.  Here $N(\rho)$ is the (prime-to-$p$) Artin conductor of $\rho$.

\begin{theorem}
\label{thm:lambda}
Fix a quadratic imaginary field $K/\Q$ and a modular residual representation  $\rhob$ such that $(\rhob,N(\rhob)^-)$ satisfies $\ram$.  Let $S(\rhob)$ denote the collection of newforms $f$ such that $\rhob_f \cong \rhob$ and  $N(\rho_f)/N(\rhob)$ is squarefree and divisible only by primes that are split in $K$.
\begin{enumerate}
\item The value of $\lambda(\sf)$ is the same for all $f \in S(\rhob)$ such that $N(\rho_f) = N(\rhob)$; denote this common value by $\lambda(\rhob)$.
\item For arbitrary $f \in S(\rhob)$, we have
$$
\lambda(\sf) = \lambda(\rhob)  +  \!\!\!\! \sum_{\ell \mid \frac{N(\rho_{f})}{N(\rhob)}} \delta_\ell(f)
$$
where $\delta_\ell(f)$ is a non-negative constant that only depends upon 
$\rhob_f$ and the restriction of $\rho_{f}$ to an inertia group at $\ell$.
\end{enumerate}
\end{theorem}

\begin{proof}
The structure of the local cohomology group at a prime $v$ that is finitely split in $\Ki/K$ is identical to the cyclotomic case.  For this reason, the arguments of \cite[Section 4]{EPW} go through verbatim.
\end{proof}

\vspace{.2cm}

\begin{remark}~
\begin{enumerate}
\item The constant $\delta_\ell(f)$ is explicitly described in \cite[p.\ 570]{EPW}.

\item As a consequence of the above theorem, the modular forms in $S(\rhob)$ with the smallest $\lambda$-invariant are the ones that are minimally ramified.  As one raises the level at split primes, 
the $\lambda$-invariant will increase or remain the same.

\item If one enlarges $S(\rhob)$ and allows for primes that are inert in $K$, the situation becomes dramatically different as is explained in the following theorem.
\end{enumerate}
\end{remark}

\begin{theorem}[Bertolini-Darmon]
Let $K$ be a quadratic imaginary field and let $f$ be a weight two newform of squarefree level $N=N^+N^-$ (with respect to $K$).  If $\eta_f(N^+\!,N^-)$ is a unit, then there exists a weight two newform $g$ of squarefree level $M$ (divisible by $N$) such that:
\begin{enumerate}
\item $\rhob_g \cong \rhob_f$;
\item $M/N$ is divisible by an even number of primes all of which are inert in $K$;
\item $\sg=0$.
\end{enumerate}
\end{theorem}

\begin{proof}
To prove this, one applies the arguments of \cite[Theorem 4.4]{BD} with $\varphi$ being reduction modulo $\pi$.  In the notation of \cite{BD}, $t_f$ then equals $\lambda(\lf)$ as we know that $\mu(\lf)=0$.  Their induction argument then produces a modular form $g$ satisfying (1) and (2) such that either $\lambda(\lg)=0$ or $\sg[\pi]=0$.  In either case, we have $\lambda(\sg)=0$ and since $\sg$ has no finite index submodules, we must have that $\sg=0$.

The hypothesis that $\cng{f}{N^+\!}{N^-}$ is a $p$-adic unit is used to ensure that the mod $p^n$ forms constructed from level-raising actually lift to true modular forms.
\end{proof}

\vspace{.1cm}

\begin{remark}~
The hypothesis that $\cng{f}{N^+\!}{N^-}$ is a $p$-adic unit is necessary in the above theorem.  Indeed, consider a newform $f$ such that $t_f(\ell) > 0$ and $\delta_\ell(f) > 0$ for some $\ell \mid N^+$.  Let $g$ be a newform of level $M$ such that $N \mid M$ and such that $\rhobar_g \cong \rhobar_f$.  Thus, $\rhobar_g$ is unramified at $\ell$ and, by level-lowering, there exists a form $h$ of level $M/\ell$ congruent to $g$.  Then, by Theorem \ref{thm:lambda}, we have 
$$
\lambda(\sg) = \lambda(\sh) + \delta_\ell(h).
$$   
As $\delta_\ell(f) > 0$ and $\rhobar_f \cong \rhobar_h$, we also have $\delta_\ell(h) > 0$.  In particular, $\lambda(\sg) > 0$ and thus $\sg \neq 0$.
\end{remark}

\begin{appendix}

\section{Surjectivity of global-to-local maps}

Let $K$ be a number field and let $F$ be a finite extension of $\Qp$ with ring of integers $\O$.  Let $V$ be a $F$-representation space for $G_K$ of dimension $d$ which is ramified at only finitely many primes.  Let $T$ be an 
$G_K$-stable lattice of $V$ and
set $A = V/T \cong \left( K / \O \right)^d$.  Further, fix  a finite set of places $\Sigma$ of $K$ containing all places over $p$ and $\infty$ along with all of the ramified primes for $V$.

Let $L/K$ be a (possibly infinite) Galois extension and let $\Sigma_L$ denote all of the places of $L$ sitting over a place in $\Sigma.$  For  $w \in \Sigma_L$, fix a subspace $\L_w \subseteq H^1(L_w,A)$ such that $\sigma \L_w = \L_{\sigma w}$ for $\sigma \in \Gal(L/K)$.  We refer to this as a {\it Selmer structure} for $A$ over $L$.  This Selmer structure induces a Selmer group
$$
\sel(L,A) = \ker \left( H^1(K_\Sigma/L,A) \to \prod_{w \in \Sigma_L} H^1(L_w,A)/\L_w \right).
$$

We note that if $K \subseteq M \subseteq L$ is a tower of Galois extensions, then a Selmer structure for $A$ over $L$ naturally induces one over $M$.  Indeed, let $v$ be a place of $M$ sitting over some place of $\Sigma$.  Then restriction yields a map
$$
H^1(M_v,A) \stackrel{\res}{\maps} H^1(L_w,A)^{\Gal(L_w/M_v)}
$$
for $w$ any place of $L$ over $v$ and we set
$$
\L_v := \res^{-1} \left( \L_w^{\Gal(L_w/M_v)} \right).
$$
This subspace is independent of the choice of $w$.

Let $\Ki/K$ be a $\Zp$-extension
and assume that we have a Selmer structure for $A$ over $\Ki$.  Then for each $n \geq 0$, we have the induced Selmer structure over $K_n$ and thus a Selmer group $\sel(K_n,A)$.  Moreover, directly from the definitions, we see that
$$
\sel(\Ki,A) = \varinjlim_n \, \sel(K_n,A).
$$

\begin{remark}
It was explained to us by Ralph Greenberg that the induced Selmer structure in the supersingular case ``regularizes" Kobayashi's plus/minus local conditions at $p$. 

Indeed, if $E/\Q$ is an elliptic curve with $a_p=0$,  Kobayashi's plus/minus local condition is a subspace of $H^1(\Q_{n,p},E[p^\infty])$ with $\Zp$-corank equal to a polynomial in $p$ of degree either $p^n$ or $p^{n-1}$ depending on the parity of $n$.  The value of the $\Zp$-corank is always strictly between $p^{n-2}$ and $p^n$.  
At the infinite level, the plus/minus local condition is defined as the direct limit of these finite-level local conditions and is a cofree submodule of $H^1(\Q_{\infty,p},E[p^\infty])$ of corank $1$.  

The induced Selmer structure is then obtained by taking invariants 
from the Selmer structure over $\Q_{\infty,p}$
down to $\Q_{n,p}$.   Since the local condition at the infinite level is cofree, the induced local condition at level $n$ has $\Zp$-corank $p^n$ (and is thus regularly behaved).   This regular behavior is used crucial in the theorem below. 

By comparing coranks, we see that Kobayashi's plus/minus local condition must be different from the induced local condition.  Fortunately, even if the corresponding finite-level Selmer groups differ, the limit of these Selmer groups both yield the same Selmer group over $\Qinf$.
\end{remark}

For $v$ a prime of $K$, let $\sigma_{m,v}$ denote the set of places of $K_{m}$ lying over $v$ and set
$$
\H_{v}(\Ki,A) = \underset{m}{\varinjlim} \prod_{w \in \sigma_{m,\ell}} H^{1}(K_{m,w},A)/\L_w.
$$
For any number field $L/K$, let 
$$
\delta(L,V) = \sum_{v~\text{complex}} d + \sum_{v~\text{real}} d_v^-
$$
where $v$ runs over archimedean places of $L$ and $d_v^-$ is the dimension of the $-1$ eigenspace of a complex conjugation over $v$ acting on $V$.  

If $\p$ is a prime of $K$, the $\Lambda$-corank of
$\L_{\mathfrak P}$ is independent of the choice of prime
$\mathfrak P$ of $\Ki$ lying over $\p$; we denote it by
$r_{\p}$.

\begin{proposition}
\label{prop:surj}
Assume that:
\begin{enumerate}
\item \label{1} no place of $K$ lying over $p$ splits completely in $\Ki$;
\item \label{2} $\sel(\Ki,A)$ is $\Lambda$-cotorsion;
\item \label{3} $H^0(\Ki,A^*)$ is finite where $A^* = \Hom(T,\mu_{p^\infty})$;
\item \label{4}
$$
\sum_{\p | p} r_\p = [K:\Q]d - \delta(K,V).
$$
\end{enumerate}
Then the global-to-local map 
$$
H^1(\Ki,A) \stackrel{\gamma}{\maps} \prod_{v \in \Sigma} \H_v(\Ki,A)
$$
is surjective.
\end{proposition}

Versions of this theorem appear in \cite{GV,Greenberg,TW}.
Our theorem differs from these results in that it allows for more general local conditions at $p$ (not just an ordinary condition) and also allows for the possibility of primes splitting infinitely in the $\Zp$-extension $\Ki/K$.  However, the basic structure of our argument is identical to all of these proofs.
We most closely follow \cite[Prop 2.1]{GV}, recalling their argument below and making the necessary changes as they arise.

\begin{proof}

We first note that it suffices to show that $\coker(\gamma)$ is finite as the target of $\gamma$ contains no proper finite index submodules.  Indeed, to see this, 
first consider a place $v$ of $K$ that is finitely decomposed in $\Ki$.  For $w$ a place of $\Ki$ over $v$, we have $\Gal(K_{\infty,w}/K_v) \cong \Zp$ which has cohomological dimension $1$.  Hence, $H^1(K_{\infty,w}/K_v,A)$ is divisible and thus $\H_v$ is also divisible as it is a direct sum of quotients of such groups.
For $v$ infinitely decomposed, as in Lemma \ref{l:hv}, we have
$\H_v \cong H^1(K_v,A) \otimes \Lambda^\vee$
which has no proper finite index submodules.

A strategy for showing that $\coker(\gamma)$ is finite is to show that the corresponding cokernel at level $n$ is finite with size bounded independent of $n$.  However, these finite level cokernels could be infinite if the characteristic power series of $\sel(\Ki,A)^\vee$ has $p$-cyclotomic zeroes.  To avoid this problem, we use Greenberg's trick of twisting the Galois module structure.

Namely, let $\kappa : \Gal(\Ki/K) \cong 1+p\Zp$ be an isomorphism and consider the twisted module $A_t := A \otimes \kappa^t$ for $t \in \Z$.
Since $A$ and $A_t$ are isomorphic as $G_{\Ki}$-modules, the Selmer structure for $A$ over $\Ki$ induces a Selmer structure on $A_t$ over $\Ki$ which we write as $\L_{w,t}$. The corresponding Selmer groups differ only by a twist; as $\Lambda$-modules, we have
$$
\sel(\Ki,A)(\kappa^t) \cong \sel(\Ki,A_t)
$$
where the $\Lambda$-module structure of the left hand side is twisted by $\kappa^t$.  

The Selmer structure for $A_t$ over $\Ki$ induces one over $K_n$ for each $n \geq 0$.  Let $\gamma_t$ denote the global-to-local map defining $\sel(\Ki,A_t)$ and let $\gamma_{n,t}$ denote the corresponding map defining $\sel(K_n,A_t)$.  To prove the theorem, it suffices to show that there is some $t$ such that $\coker(\gamma_{n,t})$ is finite for all $n \geq 0$ and of size bounded independent of $n$.   

We have
\begin{eqnarray}
\label{eqn:sn}
0 \maps \sel(K_n,A_t) \maps H^1(K_\Sigma/K_n,A_t) \stackrel{\gamma_{n,t}}{\maps} \prod_{v \in \Sigma_n} H^1(K_{n,v},A_t)/ \L_{v,t}
\end{eqnarray}
where $\Sigma_n$ is the set of places of $K_n$ over places of $\Sigma$.  We will now analyze the $\O$-corank of each term in this sequence.

As is argued in \cite{GV},
$$
\corank_\O H^1(K_\Sigma/K_n,A_t) \geq \delta(K_n,V)
$$
and, for all but finitely many $t$,
$$
\corank_\O \sel(K_n,A_t) = \corank_\O \prod_{\substack{v \in \Sigma_n  \\ v \nmid p}} H^1(K_{n,v},A_t)/ \L_{v,t} = 0.
$$
We note that the computation of the $\O$-corank of $\sel(K_n,A_t)$ uses the second hypothesis of the theorem.

Now consider a prime $v$ of $K_n$ sitting over $p$.  The local Euler characteristic of $A_t$ over $K_{n,v}$ is $d [K_{n,v}:\Qp] $.  For all but finitely many $t$, $H^0(K_{n,v},A_t)$ and $H^2(K_{n,v},A_t)$ are finite and thus, the $\O$-corank of $H^1(K_{n,v},A_t)$ equals $d[K_{n,v}:\Qp]$ for these $t$.  

Let $\P$ be any place of $\Ki$ over $v$.  By definition $\L_{v,t}:= \L_{\P,t}^{\Gal(K_{\infty,\P}/K_{n,v})}$ and thus has $\O$-corank at least $r_\p [K_{n,v}:K_\p]$ and exactly this value for all but finitely many $t$.  Summing over $v \mid p$ for such values of $t$ yields
\begin{align*}
\sum_{v | p} \corank_\O H^1(K_{n,v},A_t)/\L_{v,t} &= \sum_{v | p} d[K_{n,v}:\Qp] - r_\p [K_{n,v}:K_\p] \\
&= d [K_n:\Q] - \left(\sum_{\p | p} r_\p \right) [K_n:K] \\
&= \left(d[K:\Q]-\sum_{\p | p} r_\p \right)p^n \\
&= \delta(K,V) p^n \\
&= \delta(K_n,V).
\end{align*}
The second to last equality follows from the fourth assumption of the theorem and the final equality follows from a direct computation as $p \neq 2$. Combining the computations of the $\O$-corank of the terms of (\ref{eqn:sn}), we see that for all but finitely many $t$, the cokernel of $\gamma_{n,t}$ is finite.  

In \cite{GV}, using Poitou-Tate duality, it is moreover shown that 
if $\coker(\gamma_{n,t})$ is finite, then
$$
\left|\coker(\gamma_{n,t}) \right| \leq  \left| H^0(\Ki,A^*) \right|.
$$
This bound is independent of $n$ by the third hypothesis of the theorem.  Thus, $\coker(\gamma_t) = \coker(\gamma)$ is finite, proving the theorem.
\end{proof}

\end{appendix}


\begin{thebibliography}{99}

\bibitem{AU}
A.\ Abbes, E.\ Ullmo, 
{\it A\` propos de la conjecture de Manin pour les courbes elliptiques modulaires}, Compositio Math.\ {\bf 103} (1996), no.\ 3, 269--286.

\bibitem{ARS}
A.\ Agashe, K.\ Ribet, W.\ Stein,
{\it The Modular Degree, Congruence Primes and Multiplicity One}, preprint.

\bibitem{BD96}
M.\ Bertolini and H.\ Darmon,  
{\it Heegner points on Mumford-Tate curves}, 
Inventiones Math.\ {\bf 126} (1996), 413--456.

\bibitem{BDCD}
\bysame, {\it $p$-adic $L$-functions, and the Cerednik-Drinfeld uniformization},
Inventiones Math.\ {\bf 131} (1998), 453--491.

\bibitem{BD}
\bysame,
{\it  
Iwasawa's main conjecture for elliptic curves over anticyclotomic $\Z\sb p$-extensions},
Ann.\ of Math.\ (2) {\bf 162} (2005), 1--64.


\bibitem{Carayol}
H.\ Carayol, {\it Formes modulaires et repr\'esentations galoisiennes \`a 
valeurs dans un anneau local complet}, $p$-adic monodromy and the Birch
and Swinnerton--Dyer conjecture (Boston, 1991), Amer.\ Math.\ Soc., Providence,
RI, 213--237.

\bibitem{CK}
A.\ Cojocaru, E.\ Kani, 
{\it The modular degree and the congruence number of a weight 2 cusp form} Acta Arith.\ {\bf 114} (2004), no.\ 2, 159--167. 

\bibitem{cornut}
C.\ Cornut, {\it Mazur's conjecture on higher Heegner points},
Inventiones Math.\ {\bf 148} (2002), 495--523.

\bibitem{DDT}
H.\ Darmon, F.\ Diamond and R.\ Taylor,
{\it Fermat's Last Theorem}, 
Current Developments in Mathematics 1, 1995, International Press, 1--157. 

\bibitem{DI}
H.\ Darmon and A.\ Iovita,
{\it The anticyclotomic main conjecture for supersingular elliptic curves}, preprint. 

\bibitem{EPW}
M.\ Emerton, R.\ Pollack and T.\ Weston, {\it Variation of Iwasawa invariants in Hida families}, Invent.\ Math.\ {\bf 163} (2006), 523--580.

\bibitem{Finis}
T.\ Finis, {\it The $\mu$-invariant of anticyclotomic $L$-functions of
imaginary quadratic fields}, J. Reine und Angew. Math. {\bf 596} (2006), 131-152.

\bibitem{Greenberg}
R.\ Greenberg, 
{\it Iwasawa theory for elliptic curves}, 
Lecture Notes in Math. {\bf 1716} (1999), 51--144.

\bibitem{Gord}
\bysame,
{\it Iwasawa theory for $p$-adic representations}, 
Advanced Studies in Pure Math. {\bf 17} (1989), 97-137.

\bibitem{GV}
R.\ Greenberg and V.\ Vatsal, {\it On the Iwasawa invariants of elliptic curves}, Invent.\ Math.\ {\bf 142} (2000), 17--63.

\bibitem{Gross}
B.\ Gross, {\it Heights and the special values of $L$-series}, Numbr Theory (Montreal ,1985), CMS Conf.\ Proc.\ {\bf 7}, 115--187.

\bibitem{Groth}
A.\ Grothendieck (with M.\ Raynaud and D.S.\ Rim), {\it Groupes de monodromie en g\'eom\'etrie alg\'ebrique. I,} S\'eminaire de G\'eom\'etrie Alg\'ebrique du Bois-Marie 1967--1969 (SGA 7 I).
Lecture Notes in Math.\ {\bf 288},
Springer-Verlag, Berlin--New York, 1972.

\bibitem{Helm}
D.\ Helm, \emph{On maps between modular Jacobians and Jacobians of Shimura curves}, to appear in Israel J.\ Math.

\bibitem{Hida}
H.\ Hida, \emph{Modules of congruence of Hecke algebras and $L$-functions
associated with cusp forms}, Amer.\ J.\ Math.\ {\bf 110} (1988), 323--382.

\bibitem{IP} 
A.\ Iovita and R.\ Pollack,
{\it Iwasawa theory of elliptic curves at supersingular primes over $\Zp$-extensions of number fields},
to appear in J.\ Reine Angew.\ Math.

\bibitem{Khare}
C.\ Khare, \emph{On isomorphisms between deformation rings and Hecke rings},
Invent.\ Math.\ {\bf 154} (2003), 199--222.

\bibitem{Kimparity}
B.-D.\ Kim, \emph{The parity conjecture of elliptic curves at primes with supersingular reduction}, to appear in Compositio Mathematica.

\bibitem{Kimfinite}
B.-D.\ Kim, in preparation.

\bibitem{Kobayashi}
S.\ Kobayashi, \emph{Iwasawa theory for elliptic curves at supersingular
primes}, Invent.\ Math.\ {\bf 152} (2003), 1--36.

\bibitem{Kohel}
D.\ Kohel, \emph{Hecke module structure of quaternions},
Class field theory---its centenary and prospect (Tokyo, 1998),
Adv.\ Stud.\ Pure Math.\ {\bf 30}
Math.\ Soc.\ Japan, Tokyo, 2001, 171--195.

\bibitem{P} 
R.\ Pollack, 
{\it On the $p$-adic $L$-function of a modular form at a supersingular prime}, 
Duke Math.\ J., {\bf 118} (2003), 523--558.

\bibitem{Ribet}
K.\ Ribet, \emph{On modular representations of $\Gal(\Qbar/\Q)$ arising
from modular forms}, Invent.\ Math.\ {\bf 100} (1990), 431--476.

\bibitem{RT} 
K.\ Ribet and S.\ Takahashi,
{\it Parameterizations of elliptic curves by Shimura curves and by classical modular curves},  Elliptic curves and modular forms (Washington, DC, 1996) Proc.\ Nat.\ Acad.\ Sci.\ U.S.A. {\bf 94} (1997), 11110--11114. 

\bibitem{T}
S.\ Takahashi,
{\it Degrees of Parameterizations of Elliptic Curves by Shimura Curves}, J.\ of Number Theory {\bf 90} (2001), 74--88.

\bibitem{Vcong}
V.\ Vatsal, \emph{Canonical periods and congruence formulae},
Duke Math.\ J.\ {\bf 98} (1999), 397--419.

\bibitem{Vnonzero} 
\bysame, 
{\it Uniform distribution of Heegner points}, Invent.\ Math.\ {\bf 148} (2002), 1--46.

\bibitem{Vmu} 
\bysame, 
{\it Special values of anticyclotomic $L$-functions}, Duke Math.\ J.\ {\bf 116} (2003) 219--261.

\bibitem{TW}
T.\ Weston, 
{\it Iwasawa invariants of Galois deformations},
Manuscripta Math.\ {\bf 118} (2005), 161--180.

\bibitem{Wiles}
A.\ Wiles, \emph{Modular elliptic curves and Fermat's last theorem},
Ann.\ of Math.\ {\bf 141} (1995), 443--551.

\end{thebibliography}
\end{document}